\RequirePackage{fix-cm}
\documentclass[twocolumn]{svjour3}          
\smartqed  
\usepackage{graphicx}
\usepackage{epstopdf}
\usepackage{amssymb}
\usepackage{subfigure}
\usepackage{flushend,cuted}

\usepackage{amsmath}
\usepackage{url}
\usepackage[colorlinks=true,citecolor=blue,linkcolor=blue]{hyperref}
\usepackage{multirow}
\usepackage{booktabs}
\usepackage{mathrsfs}
\allowdisplaybreaks
\begin{document}

\title{An efficient differential quadrature method for fractional advection-diffusion equation}

\titlerunning{An efficient DQ method for fractional ADEs}        

\author{X. G. Zhu  \and Y. F. Nie \and W. W. Zhang}

\authorrunning{X. G. Zhu et al.} 

\institute{  X. G. Zhu  \and Y. F. Nie \and W. W. Zhang\at
               Department of Applied Mathematics, Northwestern  Polytechnical University, Xi'an, Shaanxi 710129, P. R. China \\
              \email{yfnie@nwpu.edu.cn}     
              \and
              X. G. Zhu \at
              \email{zhuxg590@yeah.net}     
              \and
              W. W. Zhang \at
              \email{wwzhang@nwpu.edu.cn}      
}

\date{Received: date / Accepted: date}

\maketitle

\begin{abstract}
This article  studies a direct numerical approach for fractional advection-diffusion equations (ADEs).
Using a set of cubic trigonometric B-splines as test functions, a differential quadrature (DQ) method is firstly
proposed for the 1D and 2D time-fractional ADEs of order $(0,1]$. The weighted
coefficients are determined, and with them, the original equation is transformed into a group of general
ordinary differential equations (ODEs), which are discretized by an effective difference scheme or
Runge-Kutta method.  The stability is investigated under a mild theoretical condition.
Secondly, based on a set of cubic B-splines, we develop a new Crank-Nicolson type DQ method for the 2D space-fractional ADEs without advection.
The DQ approximations to fractional derivatives are introduced and the values of the fractional derivatives
of B-splines are computed by deriving explicit formulas. The presented DQ methods are evaluated on five benchmark
problems and the concrete simulations of the unsteady propagation of solitons and Gaussian pulse.
In comparison with the existing algorithms in the open literature, numerical results finally illustrate the validity and accuracy.
\end{abstract}

\section{Introduction}\label{intro}
The differential equations with a fractional derivative serve as superior models in subjects
as diverse as astrophysics, chaotic dynamics, fractal network, signal processing, continuum mechanics,
turbulent flow and wave propagation \cite{ref02,ref01,19,20}.
This type of equations admit the non-local memory effects in mathematical mechanism, thereby filling in
a big gap that the classical models can not work well for some of the natural phenomena like anomalous
transport. In general, the exact solutions can seldom be represented as closed-form expressions by using elementary functions
that presents a tough challenge to derive a sufficiently valid  method concerned with analytic approximations,
so a keen interest has been attracted to design robust algorithms to investigate them in numerical perspectives.

In this article, we aim to construct an efficient method to numerically solve the general problems:
\begin{itemize}
  \item[(I)] 1D time-fractional ADEs
  \begin{align}
    \frac{\partial^\alpha u(x,t)}{\partial t^\alpha}+\kappa\frac{\partial u(x,t)}{\partial x}
        &-\varepsilon\frac{\partial^2 u(x,t)}{\partial x^2} =f(x,t),   \label{eq01}
\end{align}
with $0<\alpha\leq1$, $\kappa, \varepsilon \geq 0$, $a\leq x\leq b$, $t>0$, and the initial and boundary conditions
\begin{align}
  & u(x,0)=\psi(x),\quad a\leq x\leq b, \label{eq02}\\
  & u(a,t)=g_1(t), \quad u(b,t)=g_2(t), \quad t>0; \label{eq03}
\end{align}
\item[(II)] 2D time-fractional ADEs
\begin{gather}
\begin{aligned}\label{eq04}
&\frac{\partial^\alpha u(x,y,t)}{\partial t^\alpha}+\kappa_x\frac{\partial u(x,y,t)}{\partial x}+\kappa_y\frac{\partial u(x,y,t)}{\partial y}\\
&-\varepsilon_x\frac{\partial^2 u(x,y,t)}{\partial x^2}-\varepsilon_y\frac{\partial^2 u(x,y,t)}{\partial y^2} =f(x,y,t),
\end{aligned}
\end{gather}
with $0<\alpha\leq1$, $\kappa_x, \kappa_y, \varepsilon_x,\varepsilon_y \geq 0$, $(x,y)\in\Omega$, $t>0$,
and  the initial and boundary conditions
\begin{align}
  & u(x,y,0)=\psi(x,y),\quad (x,y)\in\Omega, \label{eq05}\\
  & u(x,y,t)=g(x,y,t), \ \ (x,y)\in\partial\Omega,\quad t>0, \label{eq07}
\end{align}
where $\Omega=\{(x,y):a\leq x\leq b, c\leq y\leq d\}$ and $\partial\Omega$ denotes its boundary;
\item[(III)] 2D space-fractional ADEs without advection
\begin{gather}
\begin{aligned}\label{xz01}
&\frac{\partial u(x,y,t)}{\partial t}-\varepsilon_x\frac{\partial^{\beta_1} u(x,y,t)}{\partial x^{\beta_1}}\\
&\qquad\qquad-\varepsilon_y\frac{\partial^{\beta_2} u(x,y,t)}{\partial y^{\beta_2}} =f(x,y,t),
\end{aligned}
\end{gather}
with $1<\beta_1, \beta_2\leq2$, $\varepsilon_x,\varepsilon_y \geq 0$, $(x,y)\in\Omega$, $t>0$,
and  the initial and boundary conditions
\begin{align}
  & u(x,y,0)=\psi(x,y),\quad (x,y)\in\Omega, \label{xz02}\\
  & u(x,y,t)=0, \ \ (x,y)\in\partial\Omega,\quad t>0, \label{xz03}
\end{align}
where $\Omega$ and $\partial\Omega$ are given as above.
\end{itemize}

In Eqs. (\ref{eq01}), (\ref{eq04}), the time-fractional derivatives are defined in Caputo sense, i.e.,
\begin{align*}
 \frac{\partial^\alpha u(x,t)}{\partial t^\alpha}=\frac{1}{\Gamma(1-\alpha)}
  \int^t_0\frac{\partial u(x,\xi)}{\partial \xi}\frac{d\xi}{(t-\xi)^\alpha},
\end{align*}
while in Eq. (\ref{xz01}), the space-fractional derivatives are defined in Riemann-Liouville sense, i.e.,
\begin{align*}
   \frac{\partial^{\beta_1} u(x,y,t)}{\partial x^{\beta_1}}&=\frac{1}{\Gamma(2-\beta_1)}\frac{\partial^2 }{\partial x^2}
    \int^x_a\frac{u(\xi,y,t)d\xi}{(x-\xi)^{\beta_1-1}},\\
   \frac{\partial^{\beta_2} u(x,y,t)}{\partial y^{\beta_2}}&=\frac{1}{\Gamma(2-\beta_2)}\frac{\partial^2 }{\partial y^2}
    \int^y_c\frac{u(x,\xi,t)d\xi}{(y-\xi)^{\beta_2-1}},
\end{align*}
and $\frac{\partial^\alpha u(x,y,t)}{\partial t^\alpha}$ is an analog of $\frac{\partial^\alpha u(x,t)}{\partial t^\alpha}$,
where $\Gamma(\cdot)$ is the Gamma function. It is noted that Eqs. (\ref{eq01})-(\ref{eq03}), (\ref{eq04})-(\ref{eq07}),
and (\ref{xz01})-(\ref{xz03}) reduce into the classical 1D or 2D advection-diffusion equations if $\alpha=1$, $\beta_1=\beta_2=2$ are fixed.

In recent decades, fractional ADEs have been notable subjects of intense research. 
Except for a few analytic solutions, various numerical methods have been done for Eqs. (\ref{eq01})-(\ref{eq03})
without advection, covering implicit difference method \cite{ref09},
high-order finite element method (FEM) \cite{ref12}, Legendre wavelets and spectral Galerkin methods \cite{08,ref11},
direct discontinuous Galerkin method \cite{ref13},
quadratic spline collocation method \cite{ref22}, cubic B-spline collocation method (CBCM) \cite{ref23},
orthogonal spline collocation method \cite{13},
pseudo-spectral method \cite{11}, high-order compact difference method \cite{ref10},
implicit radial basis function (RBF) meshless method \cite{12}, nonpolynomial and polynomial spline methods \cite{14}.
In \cite{ref21,ref15,ref17,06}, the algorithms based on shifted fractional Jacobi polynomials, Sinc functions and
shifted Legendre polynomials, Haar wavelets and the third kind Chebyshev wavelets functions were
well developed via the integral operational matrix or collocation strategy for Eqs. (\ref{eq01})-(\ref{eq03})
with variable coefficients. The Gegenbauer polynomial spectral collocation method was proposed in \cite{ref18}
for the same type of equations and a Sinc-Haar collocation method can be found in \cite{ref24}.
Uddin and Haq considered  a radial basis interpolation approach \cite{ref16}.
Cui established a high-order compact exponential difference scheme \cite{02}.
Razminia et al. proposed a DQ method for time-fractional diffusion equations
by using Lagrangian interpolation polynomials as test functions \cite{15}.
Shirzadi et al. solved the 2D time-fractional ADEs with a reaction term via a local Petrov-Galerkin meshless method \cite{07}.
Gao and Sun derived two different three-point combined compact alternating direction implicit (ADI)
schemes for Eqs. (\ref{eq04})-(\ref{eq07}) \cite{16}, both of which own high accuracy.
High-dimensional space-fractional ADEs are challenging topics in whether analytic or numerical aspects due to the complexity
and huge computing burden. The application of an numerical method to Eqs. (\ref{xz01})-(\ref{xz03}) did not have large diffusion;
for the conventional algorithms, we refer the readers to \cite{wx02,wx01,wx03,wx04,wx05,wx06} and references therein.

The trigonometric B-splines are a class of piecewise-defined functions constructed from algebraic trigonometric spaces,
which have got recognition since 1964. They are preferred to the familiar polynomial B-splines since
they often yield less errors when served as basis functions in interpolation theory.
Nevertheless, using these basis splines to set up numerical
algorithms is in its infancy and the related works are limited \cite{18,17}. In this study,
a DQ method for the general ADEs is developed with its weighted coefficients
calculated based on cubic trigonometric B-spline (CTB) functions. The basis splines are slightly modified for brevity and a few advantages.
Difference schemes and Runge-Kutta Gill's method are introduced to discretize the resulting ODEs.
The condition ensuring the stability of the time-stepping DQ method is discussed and found to be  mild.
Also, we propose a new cubic B-splines based DQ method for the 2D space-fractional diffusion equations by introducing the DQ approximations
to fractional derivatives. The weights are determined by deriving explicit formulas for the fractional derivatives of B-splines
via a recursive technique of partial integration. The approaches in presence are straight forward to apply and simple in implementation;
numerical results highlight the superiority over some previous algorithms.

The remainder is organized as follows. In Section \ref{s1}, we outline some basic definitions and the
cubic spline functions that are very useful hereinafter. In Section \ref{s2}, how to determine the weighted
coefficients based on these CTB functions is studied and a time-stepping DQ method is constructed for Eqs. (\ref{eq01})-(\ref{eq03})
and Eqs. (\ref{eq04})-(\ref{eq07}). Section \ref{s3} elaborates on its stable analysis.
In Section \ref{s4}, we suggest a spline-based DQ method for Eqs. (\ref{xz01})-(\ref{xz03}) based on a set of modified cubic B-splines by
explicitly computing the values of their fractional derivatives at sampling points.
Some numerical examples are included in Section \ref{s5}, which manifest the
effectiveness of our methods. The last section devotes to a conclusion.

\section{Preliminaries}\label{s1}
Let $M, N\in \mathbb{Z}^+$ and a time-space lattice be 
\begin{align*}
\Omega_{\tau}&=\{t_n:t_n=n\tau,\ 0\leq n\leq N \},\\
\Omega_h&=\{x_i:x_i=a+ih,\ 0\leq i\leq M\},
\end{align*}
with $\tau=T/N $, $h=(b-a)/M$ on $(0,T]\times[a,b]$.
Then, some auxiliary results are introduced for preliminaries.

\subsection{Fractional derivatives and their discretizations}
Given a good enough  $f(x,t)$, the formulas
\begin{align*}
&{^C_0}D^\alpha_tf(x,t)=\frac{1}{\Gamma(m-\alpha)}
\int^t_0\frac{\partial^m f(x,\xi)}{\partial \xi^m}\frac{d\xi}{(t-\xi)^{1+\alpha-m}}, \\
&{^{RL}_0}D^\alpha_tf(x,t)=\frac{1}{\Gamma(m-\alpha)}\frac{\partial^m}{\partial t^m}
\int^t_0\frac{f(x,\xi)d\xi}{(t-\xi)^{1+\alpha-m}}, 
\end{align*}
define the $\alpha$-th Caputo and Riemann-Liouville derivatives, respectively,
where $m-1<\alpha<m$, $m\in\mathbb{Z}^+$, and particularly, in the case of $\alpha=m$, both of them degenerate into
the $m$-th integer-order derivative.

The two frequently-used fractional derivatives  are equivalent with exactness to an additive factor, i.e.,
\begin{equation}\label{eq11}
    {^C_0}D^\alpha_tf(x,t)={^{RL}_0}D^\alpha_tf(x,t)-\sum^{m-1}_{l=0}\frac{f^{(l)}(x,0)t^{l-\alpha}}{\Gamma(l+1-\alpha)};
\end{equation}
see \cite{23,ref01} for references. Utilizing ${^{RL}_0}D^\alpha_t t^l=\frac{\Gamma(l+1)t^{l-\alpha}}{\Gamma(l+1-\alpha)}$
and a proper scheme to discretize the Riemann-Liouville derivatives on the right side of Eq. (\ref{eq11}),
a difference scheme for Caputo derivative can be
\begin{align}\label{eq12}
\begin{aligned}
{^C_0}D^\alpha_tf(x,t_n)&\cong \frac{1}{\tau^\alpha}\sum_{k=0}^{n}\omega^\alpha_kf(x,t_{n-k}) \\
&\ -\frac{1}{\tau^\alpha}\sum^{m-1}_{l=0}\sum^{n}_{k=0}\frac{\omega^\alpha_k f^{(l)}(x,0)t_{n-k}^{l}}{l!}.
\end{aligned}
\end{align}
with several valid alternatives of the discrete coefficients $\{\omega^\alpha_k\}_{k=0}^{n}$  \cite{66}.
 Typically, we have
\begin{equation}\label{ez06}
   \omega^\alpha_k=(-1)^k\binom\alpha k=\frac{\Gamma{(k-\alpha)}}{\Gamma{(-\alpha)}\Gamma{(k+1)}}, \ \ k\geq0,
\end{equation}
whose truncated error is $\mathscr{R}_\tau=\mathscr{O}(\tau)$, and 
\begin{equation}\label{eq08}
\omega^\alpha_k=\bigg(\frac{11}{6}\bigg)^\alpha\sum^k_{p=0}\sum^p_{q=0}\mu^q\overline{\mu}^{p-q}l^\alpha_ql^\alpha_{p-q}l^\alpha_{k-p},\ \ k\geq0,
\end{equation}
with $\mu=\frac{4}{7+\sqrt{39}\textrm{i}}$, $\overline{\mu}=\frac{4}{7-\sqrt{39}\textrm{i}}$, $\textrm{i}=\sqrt{-1}$, and
\begin{equation}\label{ez08}
  l^\alpha_0=1,\ \  l^\alpha_k=\bigg(1-\frac{\alpha+1}{k}\bigg)l^\alpha_{k-1},\ \ k\geq1,
\end{equation}
in which case, the truncated error fulfills $\mathscr{R}_\tau=\mathscr{O}(\tau^3)$. Actually, Eq. (\ref{ez08}) is the
recursive relation of Eq. (\ref{ez06}). Moreover, the coefficients $\{\omega^\alpha_k\}_{k=0}^{n}$ in Eq. (\ref{ez06}) satisfy
\begin{itemize}
   \item[(i)] $\omega^\alpha_0=1, \quad \omega^\alpha_k< 0$, \ \ $\forall k\geq 1$,
   \item[(ii)] $\sum_{k=0}^{\infty}\omega^\alpha_k=0, \quad \sum_{k=0}^{n-1}\omega^\alpha_k>0$.
\end{itemize}
These properties are easily obtained from \cite{ref01}.

Reset $0<\alpha<1$, (\ref{eq12}) thus turns into
\begin{align}\label{eq15}
\begin{aligned}
{^C_0}D^\alpha_tf(x,t_n)&= \frac{1}{\tau^\alpha}\sum_{k=0}^{n-1}\omega^\alpha_kf(x,t_{n-k})\\
&\ -\frac{1}{\tau^\alpha}\sum_{k=0}^{n-1}\omega^\alpha_kf(x,0)+\mathscr{R}_\tau.
\end{aligned}
\end{align}
It is noteworthy that Eq. (\ref{eq15}) gives a smooth transition to the classic schemes when $\alpha=1$, for instance,
Eq. (\ref{eq15}) would be the four-point backward difference scheme if $\alpha=1$ and $\{\omega^\alpha_k\}_{k=0}^{n}$ are chosen to be the ones
in Eq. (\ref{eq08}), because these coefficients also fulfill $\sum_{k=0}^{\infty}\omega^\alpha_k=0$ and vanish
apart from $\omega^\alpha_0$, $\omega^\alpha_1$, $\omega^\alpha_2$ and $\omega^\alpha_3$.

\subsection{Cubic spline functions}
Let $x_{-i}=a-ih$, $x_{M+i}=b+ih$, $i=1,2,3$ be the six ghost knots outside $[a,b]$. Then the
desirable CTB basis functions $\{CTB_m(x)\}_{m=-1}^{M+1}$ are defined as \cite{17,21}
\begin{align*}
CT{B_m}(x) = \frac{1}{\chi }\left\{ \begin{array}{l}
\phi_1(x),\qquad x \in [{x_{m - 2}},{x_{m - 1}})\\
\phi_2(x),\qquad x \in [{x_{m - 1}},{x_m})\\
\phi_3(x), \qquad x \in [{x_m},{x_{m + 1}})\\
\phi_4(x),\qquad x \in [{x_{m + 1}},{x_{m + 2}})\\
0, \qquad\quad\quad {\rm{otherwise}}
\end{array} \right.
\end{align*}
where
\begin{align*}
&\phi_1(x)={p^3}({x_{m - 2}}),  \\
&\phi_2(x)=q({x_{m + 2}}){p ^2}({x_{m-1}})+p^2({x_{m - 2}})q ({x_m})\\
   &\qquad\quad\quad+p({x_{m-2}})p ({x_{m - 1}})q({x_{m + 1}}),\\
&\phi_3(x)=p({x_{m - 2}}){q^2}({x_{m + 1}})+q^2({x_{m + 2}})p ({x_m})\\
    &\qquad\quad\quad+ p({x_{m - 1}})q({x_{m + 1}})q({x_{m + 2}}),\\
&\phi_4(x)={q ^3}({x_{m + 2}}),
\end{align*}
with the notations
\begin{align*}
&p({x_m})=\sin\bigg({\frac{{x-{x_m}}}{2}}\bigg),\\
&q({x_m})=\sin\bigg({\frac{{{x_m}-x}}{2}}\bigg),\\
&\chi =\sin \bigg({\frac{h}{2}}\bigg)\sin(h)\sin\bigg({\frac{{3h}}{2}}\bigg).
\end{align*}
The values of $CTB_m(x)$ at each knot are given by
\begin{small}
\begin{align}\label{eq16}
CT{B_m}(x_i)=\left\{
\begin{aligned}
&\sin^2\bigg(\frac{h}{2}\bigg)\csc(h)\csc\bigg(\frac{3h}{2}\bigg),\ i=m \pm 1 \\
&\frac{2}{1+2\cos(h)},\ \ i = m \\
&0,\ \ \rm{otherwise}
\end{aligned}
\right.
\end{align}
\end{small}
and the values of $CTB'_m(x)$ at each knot are given by
\begin{align}\label{eq17}
CT{B'_m}(x_i)=\left\{
\begin{aligned}
&\frac{3}{4}\csc \bigg( {\frac{{3h}}{2}} \bigg),\ \ i=m - 1\\
&-\frac{3}{4}\csc \bigg( {\frac{{3h}}{2}} \bigg),\ \ i = m+1\\
&0.\ \ \rm{otherwise}
\end{aligned}
\right.
\end{align}

Using the same grid information, the cubic B-spline basis functions $\{B_m(x)\}_{m=-1}^{M+1}$ are defined by
\begin{align*}
{B_m}(x) = \frac{1}{h^3}\left\{ \begin{array}{l}
\varphi_1(x),\qquad x \in [{x_{m - 2}},{x_{m - 1}})\\
\varphi_2(x),\qquad x \in [{x_{m - 1}},{x_m})\\
\varphi_3(x), \qquad x \in [{x_m},{x_{m + 1}})\\
\varphi_4(x),\qquad x \in [{x_{m + 1}},{x_{m + 2}})\\
0, \qquad\quad\quad {\rm{otherwise}}
\end{array} \right.
\end{align*}
with the piecewise functions
\begin{align*}
&\varphi_1(x)=(x-x_{m - 2})^3,  \\
&\varphi_2(x)=(x-x_{m - 2})^3-4(x-x_{m - 1})^3,\\
&\varphi_3(x)=(x_{m + 2}-x)^3-4(x_{m + 1}-x)^3,\\
&\varphi_4(x)=(x_{m + 2}-x)^3.
\end{align*}

Both $\{CTB_m(x)\}_{m=-1}^{M+1}$ and $\{B_m(x)\}_{m=-1}^{M+1}$ are locally compact and twice continuously differentiable on $[a,b]$.
Since the knots $x_{-1}$, $x_{M+1}$ lie beyond $[a,b]$ and the weights in relation to the B-splines
at both ends do not participate in practical computation, hereunder, as in \cite{MJ} for cubic B-splines, we modify the CTBs by
\begin{align}\label{ref22}
\left\{\begin{aligned}
&MT{B_0}(x) = CT{B_0}(x) + 2CT{B_{ - 1}}(x),\\
&MT{B_1}(x) = CT{B_1}(x) - CT{B_{ - 1}}(x),\\
&MT{B_m}(x) = CT{B_m}(x),\;\;m = 2,3, \ldots ,M - 2,\\
&MT{B_{M - 1}}(x) = CT{B_{M - 1}}(x) - CT{B_{M + 1}}(x),\\
&MT{B_M}(x) = CT{B_M}(x) + 2CT{B_{M + 1}}(x),
\end{aligned} \right.
\end{align}
for simplicity, which will result in a strictly tri-diagonal algebraic system after discretization 
on the uniform grid $\Omega_h$. $\{MTB_m(x)\}_{m=0}^{M}$ are also linearly independent and constitute a family of basis elements of a spline space.

\section{Description of CTB based DQ method}\label{s2}
On a 2D domain $[a,b]\times[c,d]$, letting $M_x, M_y\in \mathbb{Z}^+$, we add a spatial lattice
with equally spaced grid points with spacing of $h_x=(b-a)/M_x$ in $x$-axis  and $h_y=(d-c)/M_y$
in $y$-axis, i.e.,
\begin{align*}
\Omega_x&=\{x_i:x_i=a+ih_x,\ 0\leq i\leq M_x\},\\
\Omega_y&=\{y_j:y_j=c+jh_y,\ 0\leq j\leq M_y\}.
\end{align*}
DQ method is understood as a numerical technique for finding the approximate solutions of differential equations that reduces
the original problem to those of solving a system of algebraic or ordinary differential equations via
replacing the spatial partial derivatives by the representative weighted combinations of the functional values at certain grid points on the whole domain \cite{26}.
The key procedure of such method lies in the determination of its weights and the selection of the test functions
whose derivative values are explicit at the prescribed discrete grid points. As requested, we let
\begin{align}\label{eq18}
\frac{\partial ^su(x_i,t)}{\partial x^s}
\cong\sum\limits_{j=0}^M {a_{ij}^{(s)}u(x_j,t)},\ \ 0\leq i\leq M,
\end{align}
while for 2D problems, we let
\begin{align}
&\frac{\partial ^su(x_i,y_j,t)}{\partial x^s}
 \cong\sum\limits_{m=0}^{M_x} {a_{im}^{(s)}u(x_m,y_j,t)}, \label{eq19}\\
&\frac{\partial ^su(x_i,y_j,t)}{\partial y^s}
 \cong\sum\limits_{m=0}^{M_y} {b_{jm}^{(s)}u(x_i,y_m,t)},\label{eq20}
\end{align}
where $s\in\mathbb{Z}^+$, $0\leq i\leq M_x$, $0\leq j\leq M_y$ and $a_{ij}^{(s)}$, $a_{im}^{(s)}$, $b_{jm}^{(s)}$ are the
weighted coefficients allowing us to approximate the $s$-th derivatives or partial derivatives at
the given grid points in the DQ methods.

\subsection{The calculation of weighted coefficients}
In the sequel, we apply $\{MTB_m(x)\}_{m=0}^{M}$ to calculate the 1D, 2D unknown weights.
Putting $s=1$ and substituting these basis splines into Eq. (\ref{eq18}), we get
\begin{align*}
\frac{\partial MT{B_m}(x_i)}{\partial x}=\sum^M_{j=0}a_{ij}^{(1)}MT{B_m}(x_j),\ \ 0\leq i,m\leq M,
\end{align*}
with the weighted coefficients of the first-order derivative $a_{ij}^{(1)}$, $0\leq i,j\leq M$, yet to be determined.
In view of  (\ref{ref22}) and the properties (\ref{eq16})-(\ref{eq17}),
some manipulations on the above equations yield the matrix-vector forms
\begin{align}\label{eq21}
\left\{\begin{aligned}
&\textbf{A}\textbf{a}^{(1)}_0=\textbf{Z}_0,\\
&\textbf{A}\textbf{a}^{(1)}_1=\textbf{Z}_1,\\
&\qquad\ \vdots \\
&\textbf{A}\textbf{a}^{(1)}_M=\textbf{Z}_M,\\
\end{aligned} \right.
\end{align}
where $\textbf{A}$ is the $(M+1)\times(M+1)$ coefficient matrix
\begin{align*}
\textbf{A}=\left( \begin{array}{cccccc}
A_0+2A_1& A_1 & & & &\\
0& A_0 & A_1 & & & \\
 & A_1 & A_0 & A_1 & & \\
 & &\ddots&\ddots&\ddots& \\
 & & & A_1 & A_0 & 0\\
 & & & & A_1 & A_0+2A_1
\end{array} \right),
\end{align*}
\begin{align*}
A_0=\frac{2}{1+2\cos(h)},\ \
A_1=\sin^2\bigg(\frac{h}{2}\bigg)\csc(h)\csc\bigg(\frac{3h}{2}\bigg),
\end{align*}
$\textbf{a}^{(1)}_k$, $0\leq k\leq M$, are the weighted coefficient vectors at $x_k$, i.e.,
$\textbf{a}^{(1)}_k=[a_{k0}^{(1)},a_{k1}^{(1)},\ldots,a_{kM}^{(1)}]^{\rm{T}}$,
and the right-side vectors $\textbf{Z}_k$ at $x_k$, $0\leq k\leq M$,  are as follows
\begin{align*}
&\textbf{Z}_0 = \left(\begin{array}{c}
-2z\\ 2z\\ 0\\ 0\\ \vdots \\  0\\  0
\end{array} \right),\ \
\textbf{Z}_1=\left(\begin{array}{c}
-z\\ 0\\ z\\ 0\\ \vdots \\ 0\\ 0
\end{array}\right), \cdots, \\
&\textbf{Z}_{M-1}=\left( \begin{array}{c}
0\\ 0\\ \vdots \\ 0\\ -z\\ 0\\ z
\end{array}\right),\ \
\textbf{Z}_M=\left( \begin{array}{c}
0\\ 0\\ \vdots \\ 0\\ 0\\ -2z\\ 2z
\end{array}\right),
\end{align*}
with $z=\frac{3}{4}\csc\left(\frac{3h}{2}\right)$, respectively. Thus, $a_{ij}^{(1)}$ are obtained by solving Eqs. (\ref{eq21}) for each point $x_i$.
There are two different way to derive the weighted coefficients $a_{ij}^{(2)}$ of the second-order derivative:
(i) do a similar fashion as above by putting $s=2$ in Eq. (\ref{eq18}) and solve an algebraic system for each grid point;
(ii) find the weighted coefficients $a_{ij}^{(s)}$, $s\geq2$, corresponding to the high-order derivatives in a recursive style \cite{27}, i.e.,
\begin{align*}
&a_{ij}^{(s)}=s\Bigg( {a_{ii}^{(s-1)}a_{ij}^{(1)}-\frac{a_{ij}^{(s-1)}}{x_i-x_j}}\Bigg),\ \ i\ne j,\ 0\leq i\leq M,\\
&a_{ii}^{(s)}=-\sum\limits_{j=0,j\ne i}^M {a_{ij}^{(s)}} ,\ \ i = j,
\end{align*}
which includes $s=2$ as a special case. The former would be less efficient since
the associated equations have to be solved as priority, so the latter one will be selected during our entire computing process.
Proceeding as before via replacing $\Omega_h$ by $\Omega_x$, $\Omega_y$ leads  
to a 2D generalization to get $a^{(1)}_{im}$, $b^{(1)}_{jm}$ of the first-order partial
derivatives with regard to variables $x$, $y$ in Eqs. (\ref{eq19})-(\ref{eq20}) and by them, the following relationships can further be applied, i.e.,
\begin{align*}
&a_{im}^{(s)}=s\Bigg( {a_{ii}^{(s-1)}a_{im}^{(1)}-\frac{a_{im}^{(s-1)}}{x_i-x_m}}\Bigg),\ \ i\ne m,\ 0\leq i\leq M,\\
&a_{ii}^{(s)}=-\sum\limits_{m=0,m\ne i}^{M_x}{a_{im}^{(s)}} ,\ \ i = m,\\
&b_{jm}^{(s)}=s\Bigg( {b_{jj}^{(s-1)}b_{jm}^{(1)}-\frac{b_{jm}^{(s-1)}}{y_j-y_m}}\Bigg),\ \ j\ne m,\ 0\leq j\leq M,\\
&b_{jj}^{(s)}=-\sum\limits_{m=0,m\ne j}^{M_y}{b_{jm}^{(s)}} ,\ \ j = m,
\end{align*}
to calculate $a^{(s)}_{im}$, $b^{(s)}_{jm}$ with $s\geq2$.

A point worth noticing is that $A_0,\ A_1>0$, when $0<h<1$, $0<h_x,\ h_y<1$.
Since $A_0,\ A_1$ can be deemed to be the functions of  $h$,  we obtain their derivatives
\begin{align*}
   A'_0&=\frac{4\sin(h)}{(1+2\cos(h))^2},\\
   A'_1&=\frac{\sec(\frac{h}{2})\tan(\frac{h}{2})(5+6\cos(h))}{4(1+2\cos(h))^2}.
\end{align*}
On letting $0<h<1$,  both are proved to be larger than zero, i.e.,
$A_0,\ A_1$ are the increasing functions with respect to $h$. On the other hand, there exist
$A_0(0)=0.6667$, $A_1(1)=0.2738$. Then, it suffices to show
\begin{equation*}
    \frac{2}{1+2\cos(h)}>2\sin^2\bigg(\frac{h}{2}\bigg)\csc(h)\csc\bigg(\frac{3h}{2}\bigg),
\end{equation*}
which implies $A_0>2A_1$, and thus $\textbf{A}$ is a strictly diagonally dominant tri-diagonal matrix.
Hence, Thomas algorithm can be applied to tackle the algebraic equations as Eqs. (\ref{eq21}), which
simply requires the arithmetic operation cost $\mathscr{O}(M+1)$ and would greatly economize on the memory and
computing time in practice.

\subsection{Construction of CTB based DQ method}
In this subpart, a DQ method based on $\{MTB_m(x)\}_{m=0}^{M}$ (MCTB-DQM) is constructed for Eqs. (\ref{eq01})-(\ref{eq03}) and Eqs. (\ref{eq04})-(\ref{eq07}).
Let $s=1$, $2$. The substitution of the weighted sums (\ref{eq18}), (\ref{eq19})-(\ref{eq20}) into the main equations gives
\begin{small}
\begin{align*}
\frac{\partial^\alpha u(x_i,t)}{\partial t^\alpha}+\kappa\sum\limits_{j=0}^M {a_{ij}^{(1)}u(x_j,t)}
 -\varepsilon\sum\limits_{j=0}^M {a_{ij}^{(2)}u(x_j,t)}=f(x_i,t),
\end{align*}
\end{small}
with $i=0,1,\cdots,M$, and
\begin{small}
\begin{align*}
&\frac{\partial^\alpha u(x_i,y_j,t)}{\partial t^\alpha}\!+\!\kappa_x\sum\limits_{m=0}^{M_x}{a_{im}^{(1)}u(x_m,y_j,t)}
\!+\!\kappa_y\sum\limits_{m=0}^{M_y}{b_{jm}^{(1)}u(x_i,y_m,t)}\\
&\!-\!\varepsilon_x\sum\limits_{m=0}^{M_x}{a_{im}^{(2)}u(x_m,y_j,t)}
\!-\!\varepsilon_y\sum\limits_{m=0}^{M_y}{b_{jm}^{(2)}u(x_i,y_m,t)}\!=\!f(x_i,y_j,t),
\end{align*}
\end{small}
with $i=0,1,\cdots,M_x$, $j=0,1,\cdots,M_y$, which are indeed a group of $\alpha$-th ODEs associated with the
boundary constraints (\ref{eq03}), (\ref{eq07}), and involve $\alpha\in(0,1)$ and $\alpha=1$ as two
separate cases. In what follows, we employ
\begin{align*}
&u_i^n=u(x_i,t_n),\quad u_{ij}^n=u(x_i,y_j,t_n),\\
&f^n_{i}=f(x_i,t_n),\quad f^n_{ij}=f(x_i,y_j,t_n),\\
&g_1^n=g_1(t_n),\quad g_2^n=g_2(t_n),\quad g_{ij}^n=g(x_i,y_j,t_n),
\end{align*}
for the ease of exposition, where $n=0,1,\cdots,N$.

\subsubsection{The case of fractional order}
\vspace{-2mm}
Discretizing the ODEs above by the difference scheme (\ref{eq15}) and imposing boundary constraints, we have
\begin{align}
\left\{\begin{aligned}\label{eq22}
&\omega^\alpha_0U_i^n+\kappa\tau^\alpha\sum\limits_{j=1}^{M-1}a_{ij}^{(1)}U_j^n
    -\varepsilon\tau^\alpha\sum_{j=1}^{M-1}a_{ij}^{(2)}U_j^n\\
&=-\sum_{k=1}^{n-1}\omega^\alpha_kU_i^{n-k}\!+\!\sum_{k=0}^{n-1}\omega^\alpha_kU_i^0+ \tau^\alpha G^n_i,
\end{aligned}\right.
\end{align}
with $i=1,2,\cdots,M-1$ and
\begin{equation*}
  G^n_i=f^n_i-\kappa\big(a_{i0}^{(1)}g^n_1+a_{iM}^{(1)}g^n_2\big)+\varepsilon\big(a_{i0}^{(2)}g^n_1+a_{iM}^{(2)}g^n_2\big),
\end{equation*}
for Eqs. (\ref{eq01})-(\ref{eq03}), and the following scheme
\begin{small}
\begin{align}
\left\{\begin{aligned}\label{eq23}
&\omega^\alpha_0U_{ij}^n+\kappa_x\tau^\alpha\sum_{m=1}^{M_x-1}a_{im}^{(1)}U_{mj}^n+\kappa_y\tau^\alpha\sum_{m=1}^{M_y-1}b_{jm}^{(1)}U_{im}^n\\
&\ -\varepsilon_x\tau^\alpha\sum_{m=1}^{M_x-1}a_{im}^{(2)}U_{mj}^n-\varepsilon_y\tau^\alpha\sum_{m=1}^{M_y-1}b_{jm}^{(2)}U_{im}^n\\
&\ =-\sum_{k=1}^{n-1}\omega^\alpha_kU_{ij}^{n-k}+\sum_{k=0}^{n-1}\omega^\alpha_kU_{ij}^0+ \tau^\alpha G^n_{ij},
\end{aligned}\right.
\end{align}
\end{small}
with $i=1,2,\cdots,M_x-1$, $j=1,2,\cdots,M_y-1$, and
\begin{small}
\begin{align*}
  G^n_{ij}\!=&f^n_{ij}\!-\!\kappa_x\big(a_{i0}^{(1)}g^n_{0j}+a_{iM_x}^{(1)}g^n_{M_xj}\big)
           \!-\!\kappa_y\big(b_{j0}^{(1)}g^n_{i0}+b_{jM_y}^{(1)}g^n_{iM_y}\big)\\
          &+\varepsilon_x\big(a_{i0}^{(2)}g^n_{0j}+a_{iM_x}^{(2)}g^n_{M_xj}\big)
          +\varepsilon_y\big(b_{j0}^{(2)}g^n_{i0}+b_{jM_y}^{(2)}g^n_{iM_y}\big),
\end{align*}
\end{small}
for Eqs. (\ref{eq04})-(\ref{eq07}). The Eqs. (\ref{eq22})-(\ref{eq23}) can further be rewritten in
matrix-vector forms, for instance, letting
\begin{align*}
\textbf{U}^n&=[U^n_{11},\ldots,U^n_{M_{x}-1,1},U^n_{12},\ldots,U^n_{M_{x}-1,M_{y}-1}]^T,\\
\textbf{G}^n&=[G^n_{11},\ldots,G^n_{M_{x}-1,1},G^n_{12},\ldots,G^n_{M_{x}-1,M_{y}-1}]^T,
\end{align*}
for Eqs. (\ref{eq23}), we have
\begin{equation}\label{ez81}
   \omega^\alpha_0\textbf{U}^n+\tau^\alpha\textbf{K}\textbf{U}^n=-\sum_{k=1}^{n-1}\omega^\alpha_k\textbf{U}^{n-k}
    +\sum_{k=0}^{n-1}\omega^\alpha_k\textbf{U}^0+\tau^\alpha\textbf{G}^n,
\end{equation}
where
\begin{equation*}
    \textbf{K}=\kappa_x\textbf{I}_y\otimes \textbf{W}^{1}_x+\kappa_y \textbf{W}^{1}_y\otimes\textbf{I}_x
         -\varepsilon_x\textbf{I}_y\otimes \textbf{W}^{2}_x-\varepsilon_y \textbf{W}^{2}_y\otimes\textbf{I}_x,
\end{equation*}
with $\textbf{I}_x$, $\textbf{I}_y$ being the identity matrices in $x$- and $y$-axis,
``$\otimes$'' being Kronecker product, and
\begin{align*}
\textbf{W}^{c}_{z}
=\left( \begin{array}{llll}
w^{(c)}_{11}&w^{(c)}_{12} &\cdots &w^{(c)}_{1,M_z-1}\\
w^{(c)}_{21}&w^{(c)}_{22} &\cdots &w^{(c)}_{2,M_z-1}  \\
 \vdots &\vdots &\ddots&\vdots \\
w^{(c)}_{M_z-1,1}&w^{(c)}_{M_z-1,2}&\cdots &w^{(c)}_{M_z-1,M_z-1}
\end{array} \right), \ c=1,2,
\end{align*}
in which, $z=x$ if $w=a$ while $z=y$ if $w=b$.
The initial states are got from Eqs. (\ref{eq02}), (\ref{eq05}).
As a result, the approximate solutions are obtained via performing the iteration in Eqs. (\ref{eq22})-(\ref{eq23})
until the last time level by rewriting them in matrix-vector forms first. 

\subsubsection{The case of integer order}
When $\alpha=1$, despite $\omega^1_0=1.8333$, $\omega^1_1=-3$, $\omega^1_2=1.5$, $\omega^1_3=-0.3333$, being the
coefficients of the four-point backward difference scheme, the initial values with the errors of the same convergent rate are
generally necessary to start Eqs. (\ref{eq22})-(\ref{eq23}). However, this situation would not happen if $\{\omega_k^\alpha\}_{k=0}^n$ in Eq. (\ref{ez06})
are applied. In such a case, to make the algorithm to be more cost-effective, we use Runge-Kutta Gill's method to
handle those ODEs instead, which is explicit and fourth-order convergent. Rearrange the ODEs in a unified form
\begin{equation}\label{ez09}
    \frac{\partial\textbf{u}}{\partial t}=\textbf{F}(\textbf{u}),
\end{equation}
then the DQ method is constructed as follow
\begin{small}
\begin{align}\label{eq26}
&\textbf{U}^{n}=\textbf{U}^{n-1}+\frac{1}{6}\Big[K_1+(2-\sqrt{2})K_2+(2+\sqrt{2})K_3+K_4\Big],\nonumber\\
&K_1=\tau \textbf{F}\big(t_{n-1},\textbf{U}^{n-1}\big),\nonumber\\
&K_2=\tau \textbf{F}\bigg(t_{n-1}+\frac{\tau}{2},\textbf{U}^{n-1}+\frac{K_1}{2}\bigg),\\
&K_3=\tau \textbf{F}\bigg(t_{n-1}+\frac{\tau}{2},\textbf{U}^{n-1}+\frac{\sqrt{2}-1}{2}K_1+\frac{2-\sqrt{2}}{2}K_2\bigg),\nonumber\\
&K_4=\tau \textbf{F}\bigg(t_{n-1}+\tau,\textbf{U}^{n-1}-\frac{\sqrt{2}}{2}K_2+\frac{2+\sqrt{2}}{2}K_3\bigg),\nonumber
\end{align}
\end{small}
where $\textbf{u}$, $\textbf{U}^{n}$, $n=1,2,\ldots,N$, are the unknown vectors and $\textbf{F}(\cdot)$ stands for
the matrix-vector system corresponding to the weighted sums in ODEs and contains  $a^{(s)}_{ij}$ or $a^{(s)}_{im}$,
$b^{(s)}_{jm}$, $s=1$, $2$, as its elements. Meanwhile, the boundary constraints (\ref{eq03}), (\ref{eq07}) must be imposed on
$\textbf{F}(\cdot)$ in the way as they are done for the fractional cases before we can fully run the procedures for Eqs. (\ref{eq26}).

\begin{figure}
\centering
\includegraphics[width=3.3in]{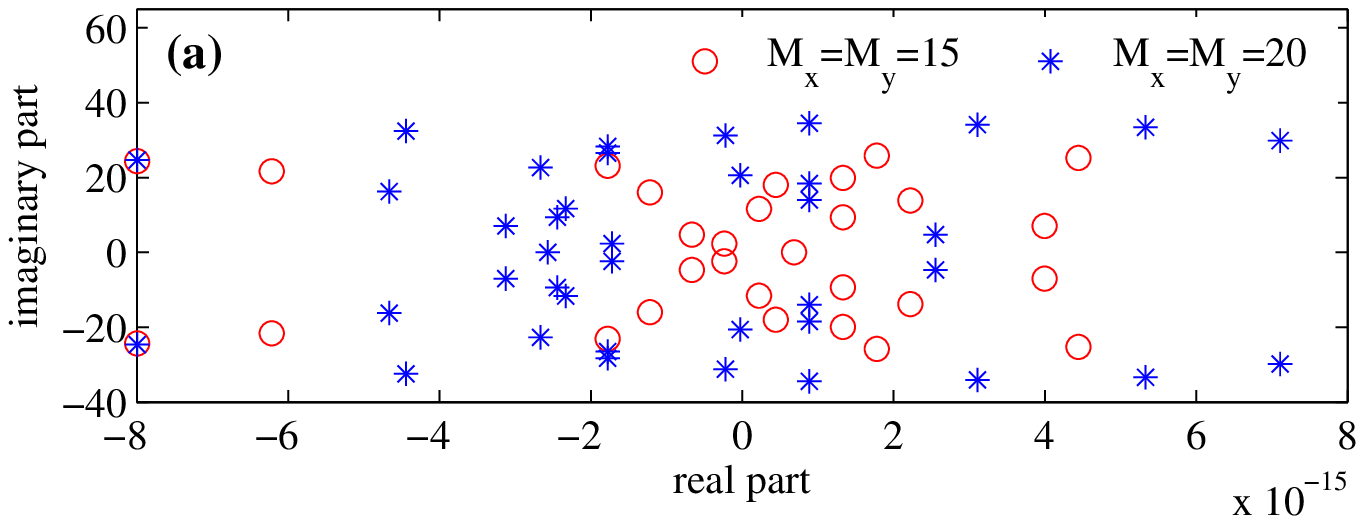}\\
\includegraphics[width=3.3in]{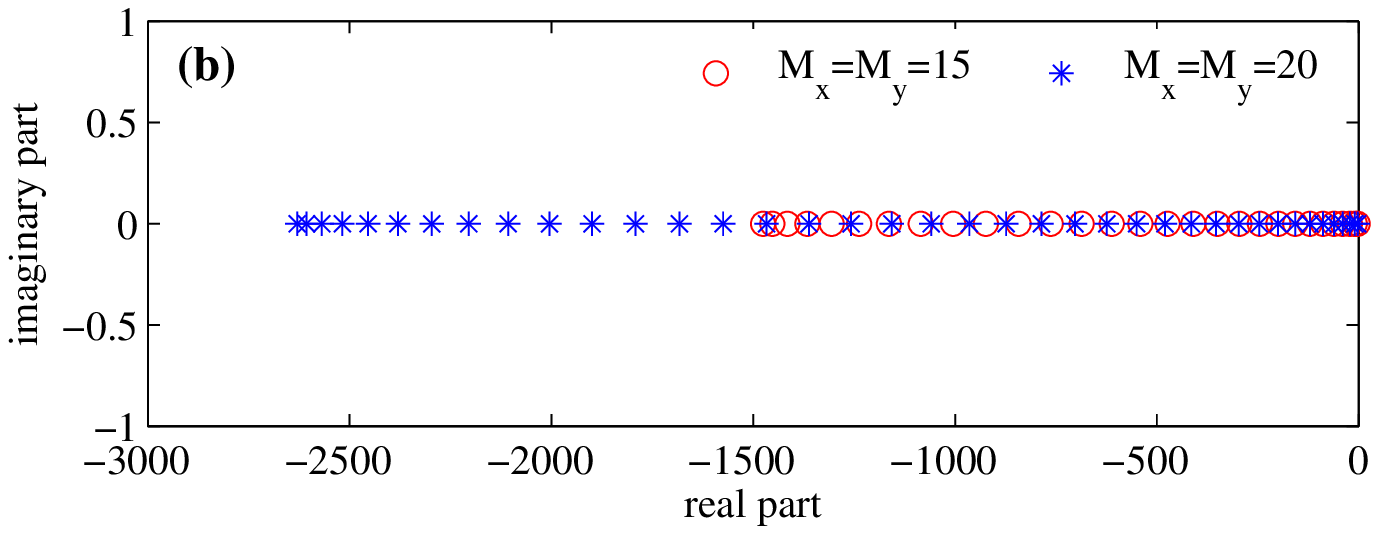}
\caption{The eigenvalues of the weighted matrices generated by DQ method when $a=c=0$, $b=d=2$: (a) $\textbf{W}^1_z$; (b) $\textbf{W}^2_z$.}\label{zfig0}
\end{figure}

\section{Stability analysis}\label{s3}
This part makes a attempt to study the matrix stability of Eqs. (\ref{ez09}) and
the numerical stability of Eqs. (\ref{eq22})-(\ref{eq23}). When $\alpha=1$,
we rewrite Eqs. (\ref{ez09}) by
\begin{equation}\label{ez12}
    \frac{\partial\textbf{u}}{\partial t}=-\textbf{K}\textbf{u} + \textbf{Q},
\end{equation}
where $\textbf{Q}$ is a vector containing the non-homogenous part and the boundary conditions, and
$\textbf{K}$ is the weighted matrix mentioned before.
 We discuss the homogeneous case. The numerical stability of an algorithm for the ODEs generated by a DQ method
relies on the stability of the ODEs themselves. Only when their solutions
are stable can a well-known method such as Runge-Kutta Gill's method
yield convergent solutions. It is enough to show their stability that the real parts of the eigenvalues of
$-\textbf{K}$ are all negative. Denote
the row vector consisting of the eigenvalues of $\textbf{W}^{c}_z$ by $\boldsymbol{\lambda}_z^c$, with $z=x,y$ and $c=1,2$. In view of
the properties of Kronecker product, the eigenvalues of $\textbf{W}^{c}_y\otimes\textbf{I}_x$, $\textbf{I}_y\otimes \textbf{W}^{c}_x$
are $\boldsymbol{\lambda}^{c}_y\otimes\textbf{e}_x$ and $\textbf{e}_y\otimes\boldsymbol{\lambda}^{c}_x$ (see \cite{xz99}), respectively,
 therefore, we have the eigenvalues of $-\textbf{K}$ in Eq. (\ref{ez12}), i.e.,
\begin{equation*}
    \boldsymbol{\lambda}=-\kappa_x\textbf{e}_y\otimes\boldsymbol{\lambda}^{1}_x-\kappa_y \boldsymbol{\lambda}^{1}_y\otimes\textbf{e}_x
         +\varepsilon_x\textbf{e}_y\otimes\boldsymbol{\lambda}^{2}_x+\varepsilon_y \boldsymbol{\lambda}^{2}_y\otimes\textbf{e}_x,
\end{equation*}
where $\textbf{e}_x$, $\textbf{e}_y$ are the row vectors of sizes $M_x+1$ and $M_y+1$, respectively, with all of their components being 1.
The exact solution of ODEs is related to $\textrm{Re}\{\boldsymbol{\lambda}\}$ and the
condition $\textrm{Re}\{\boldsymbol{\lambda}\}\leq 0$ is easy to meet because $\boldsymbol{\lambda}_z^2$ are always verified to be real and negative while
$\boldsymbol{\lambda}_z^1$ be complex with their real parts being very close to zero;
see Fig. \ref{zfig0} for example. More than that, we notice that the foregoing analysis is also valid for the 1D cases and
the phenomena appearing in Fig. \ref{zfig0} would be enhanced as the grid numbers increase.
Hence, we come to a conclusion that the ODEs are stable in most cases. 

The discussion about the numerical stability of a fully discrete DQ method is difficulty and still sparse \cite{TS0,TS1}.
In the sequel, we show the conditionally stable nature of Eqs. (\ref{eq22})-(\ref{eq23}) in the context of $L_2$-norm $||\cdot||$
and the analysis is not just applicable to the fractional case. Without loss of generality, consider the 2D cases and the
discrete coefficients $\{\omega^\alpha_k\}_{k=0}^{n}$ in Eq. (\ref{ez06}). Let $\tilde{\textbf{U}}^0$ be
the approximation of initial values $\textbf{U}^0$. Then
\begin{small}
\begin{equation}\label{ez80}
    \tilde{\textbf{U}}^n+\tau^\alpha\textbf{K}\tilde{\textbf{U}}^n=-\sum_{k=1}^{n-1}\omega^\alpha_k\tilde{\textbf{U}}^{n-k}
    +\sum_{k=0}^{n-1}\omega^\alpha_k\tilde{\textbf{U}}^0+\tau^\alpha\textbf{G}^n. 
\end{equation}
\end{small}
On subtracting Eq. (\ref{ez80}) from Eq. (\ref{ez81}) and letting $\textbf{e}^n=\textbf{U}^n-\tilde{\textbf{U}}^n$, we have
the perturbation equation
\begin{small}
\begin{equation}\label{ez82}
    \textbf{e}^n=-\sum_{k=1}^{n-1}\omega^\alpha_k(\textbf{I}+\tau^\alpha\textbf{K})^{-1}\textbf{e}^{n-k}
    +\sum_{k=0}^{n-1}\omega^\alpha_k(\textbf{I}+\tau^\alpha\textbf{K})^{-1}\textbf{e}^0,
\end{equation}
\end{small}
where $\textbf{I}$ is the identity matrix in the same size of \textbf{K}.
To prove $||\textbf{e}^n||\leq||\textbf{e}^0||$, we make the assumption
\begin{equation}\label{ez83}
    ||(\textbf{I}+\tau^\alpha\textbf{K})^{-1}||\leq 1.
\end{equation}
When $n=1$, by taking $||\cdot||$ on both sides of Eq. (\ref{ez82}), $||\textbf{e}^1||\leq ||\textbf{e}^0||$ is
trivial due to $\omega^\alpha_0=1$. Let
\begin{equation*}
    ||\textbf{e}^m||\leq ||\textbf{e}^0||, \quad m=1,2,\ldots,n-1.
\end{equation*}
Using mathematical induction, it thus follows from the properties of $\{\omega^\alpha_k\}_{k=0}^{n}$ stated in Section \ref{s1} that
\begin{small}
\begin{align*}
 ||\textbf{e}^n||&\!=\!\Bigg|\Bigg|\!-\!\sum_{k=1}^{n-1}\omega^\alpha_k(\textbf{I}\!+\!\tau^\alpha\textbf{K})^{-1}\textbf{e}^{n-k}
  +\sum_{k=0}^{n-1}\omega^\alpha_k(\textbf{I}\!+\!\tau^\alpha\textbf{K})^{-1}\textbf{e}^0\Bigg|\Bigg| \\
  &\!\leq\!\Bigg(1\!-\!\sum_{k=0}^{n-1}\omega^\alpha_k\!+\!\sum_{k=0}^{n-1}\omega^\alpha_k\Bigg)
   ||(\textbf{I}\!+\!\tau^\alpha\textbf{K})^{-1}||\max\limits_{0\leq m\leq n-1}||\textbf{e}^m|| \\
  &=||(\textbf{I}+\tau^\alpha\textbf{K})^{-1}||\max\limits_{0\leq m\leq n-1}||\textbf{e}^m|| \leq ||\textbf{e}^0||.
\end{align*}
\end{small}

\begin{figure}
\centering
\includegraphics[width=3.3in]{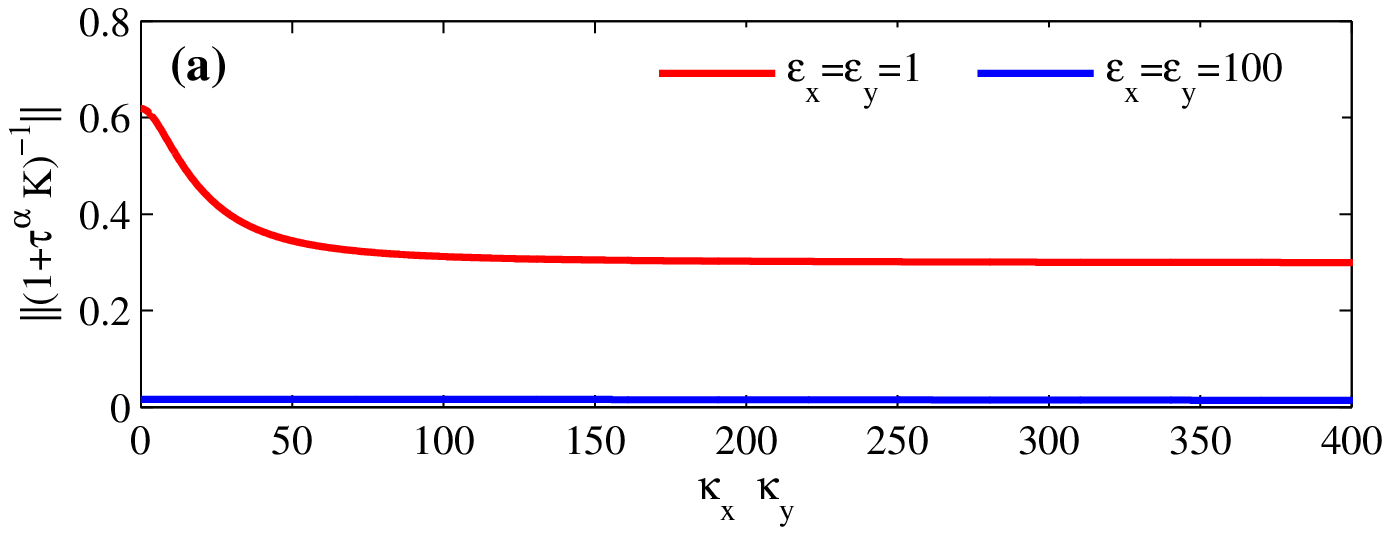}\\
\includegraphics[width=3.3in]{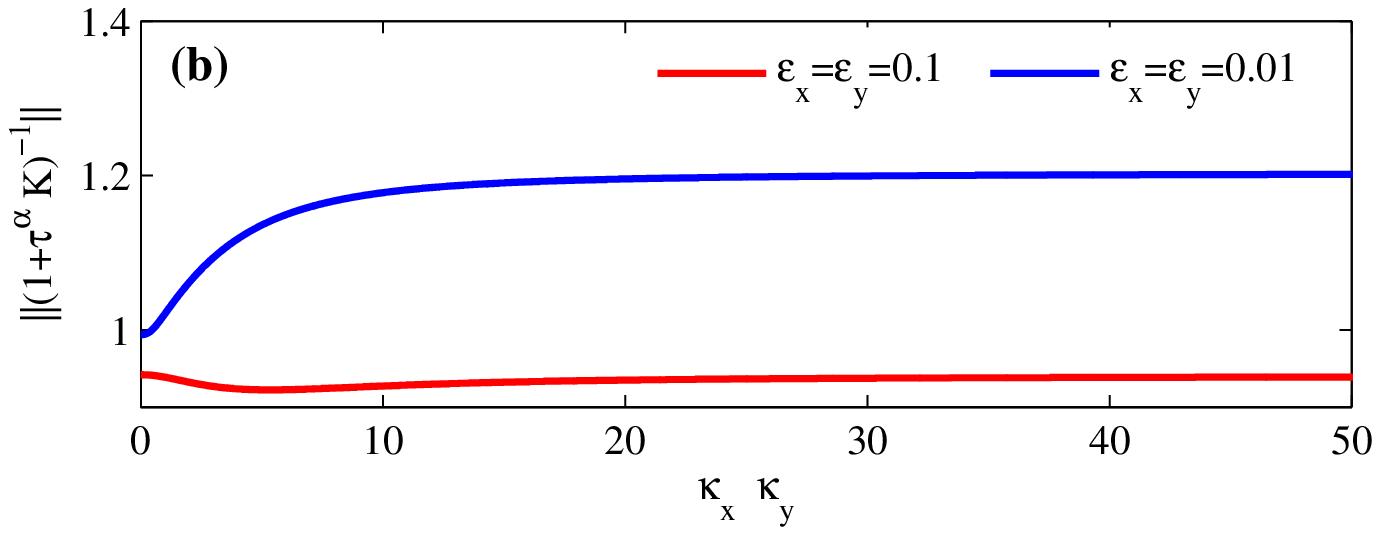}\\
\includegraphics[width=3.3in]{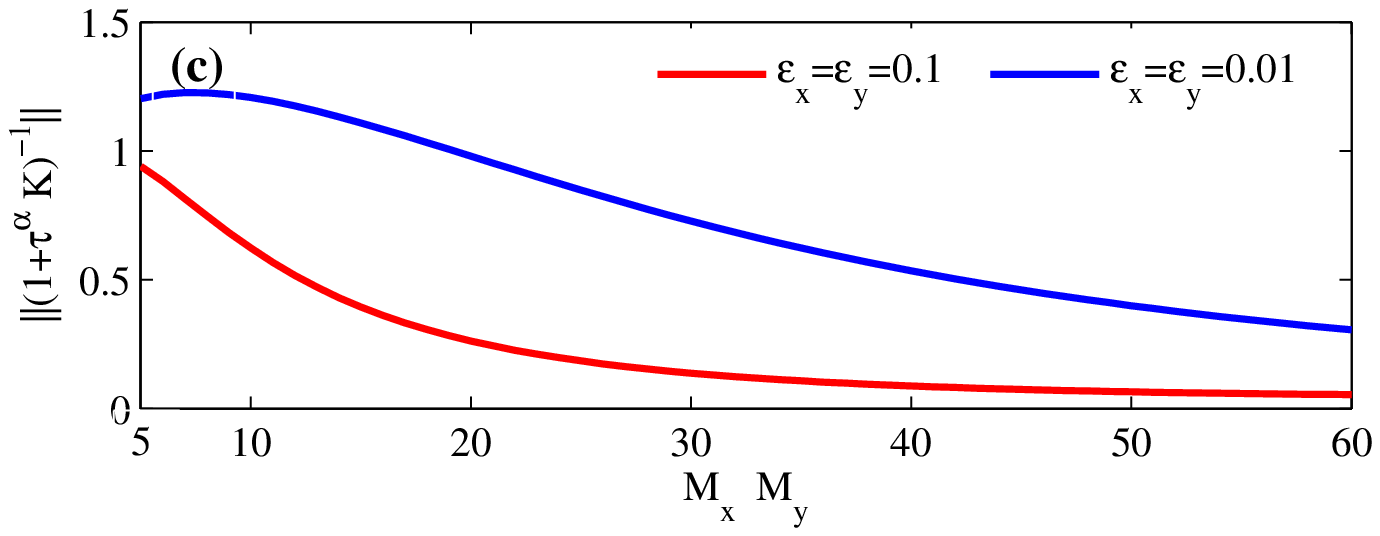}\\
\includegraphics[width=3.3in]{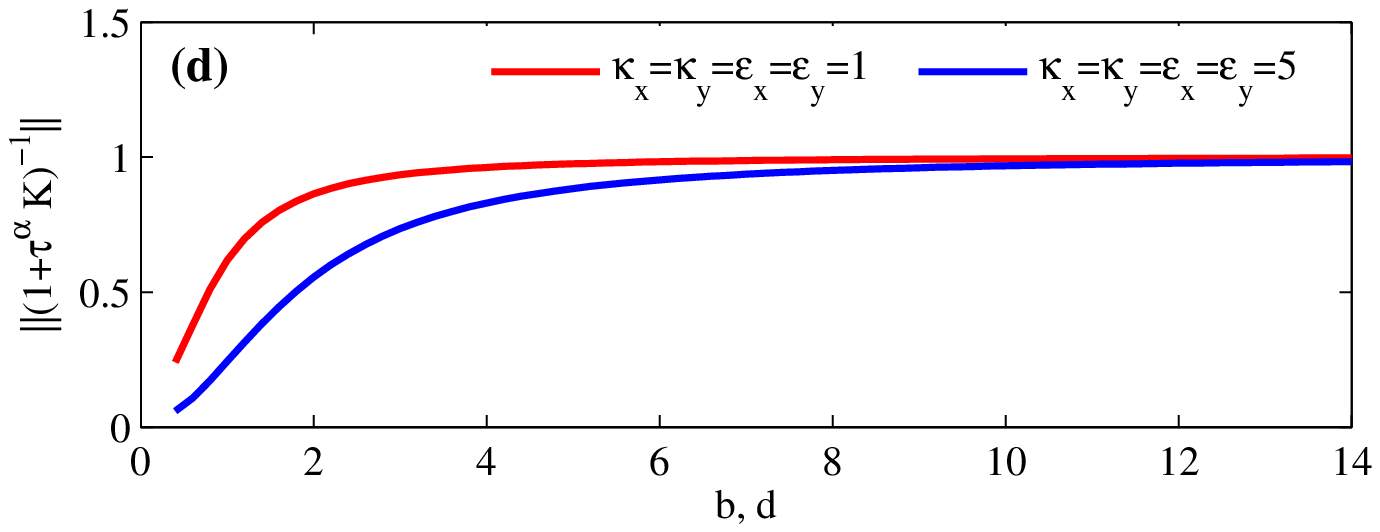}\\
\includegraphics[width=3.3in]{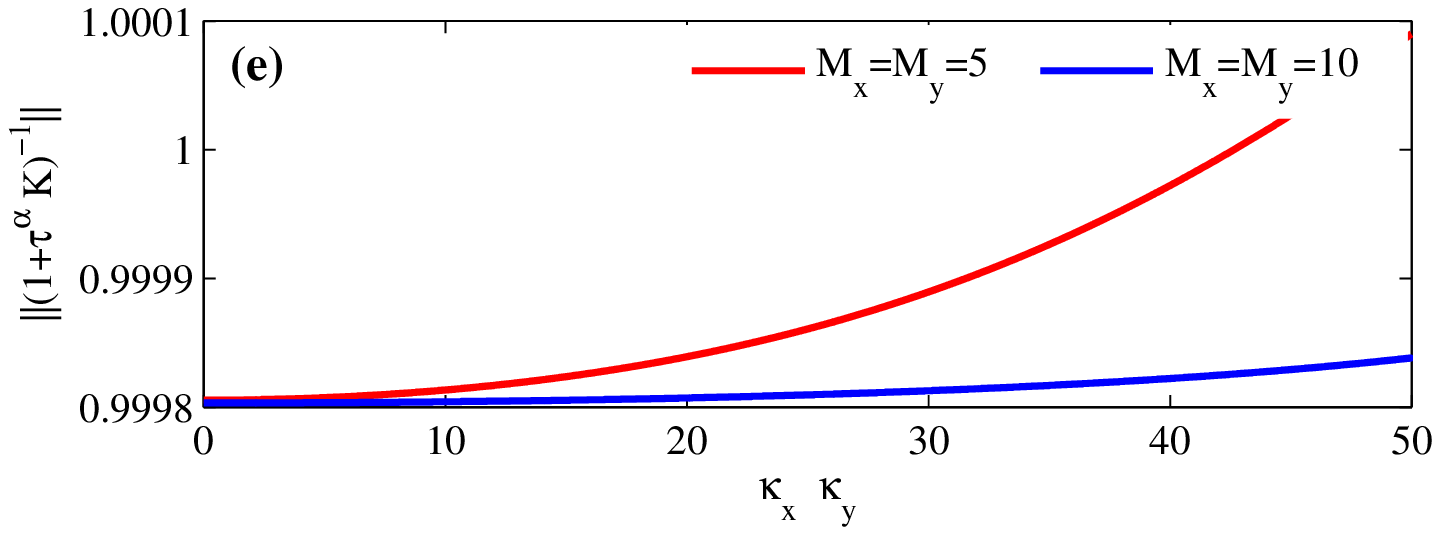}
\caption{The values of $||(\textbf{I}+\tau^\alpha\textbf{K})^{-1}||$ versus the variation of various factors:
    $\kappa_x,\kappa_y$, $\varepsilon_x,\varepsilon_y$, $M_x, M_y$, and $b, d$.}\label{zfig1}
\end{figure}

Hereinafter, we proceed with a full numerical investigation on the assumption (\ref{ez83}) to explore the
potential factors which may lead to $||(\textbf{I}+\tau^\alpha\textbf{K})^{-1}||> 1$. At first, if $\tau^\alpha$
continuously varies from 1 to 0, there holds $||(\textbf{I}+\tau^\alpha\textbf{K})^{-1}||\rightarrow 1$. However,
this process can affect the maximal ratio of the coefficients of advection and diffusivity to keep (\ref{ez83});
we leave this case to the end of the discussion.
To be more representative, we take $\tau=1.0\times10^{-3}$, $\alpha=0.5$, and $M_x=M_y=5$, unless otherwise stated. The main procedures are divided
into three steps: (i) fixing $\varepsilon_x,\varepsilon_y$, and $\Omega$, let $\kappa_x,\kappa_y$ vary and the
values of $||(\textbf{I}+\tau^\alpha\textbf{K})^{-1}||$ as the function of $\kappa_x,\kappa_y$ are plotted in (a), (b) of Fig. \ref{zfig1};
(ii) fixing $\kappa_x,\kappa_y$, and $\Omega$, let $M_x$, $M_y$ vary and the results are plotted in (c) of Fig. \ref{zfig1},
where $\kappa_x=\kappa_y=500$; (iii) fixing $\kappa_x,\kappa_y,\varepsilon_x,\varepsilon_y$, let $a=c=0$ and $b$, $d$ vary,
and the corresponding results are presented in (d) of Fig. \ref{zfig1}. It is worthy to note that $\Omega$ is the unit square
except the case of (iii), and the parameters of the same types in $x$- and $y$- axis are used as the same, for example, $\varepsilon_x=\varepsilon_y$.
Now, we consider the influence brought by $\tau$. Resetting $\tau=1.0\times10^{-10}$, let $\varepsilon_x=\varepsilon_y=1$ and $\kappa_x, \kappa_y$ vary.
The behavior of objective quantity is plotted in subfigure (e), from which, we see that the critical ratio between $\kappa_x, \kappa_y$ and
$\varepsilon_x, \varepsilon_y$ to maintain (\ref{ez83}) is about 40, far less than the case of (i), and can further be improved by increasing $M_x, M_y$.

From the foregoing discussion and figures, we summarize the conclusions as follows:
(i) if $\varepsilon_x$, $\varepsilon_y$ are not small, the tolerant ranges of $\kappa_x$, $\kappa_y$ to
    guarantee (\ref{ez83}) is quite loose and when $\varepsilon_x$, $\varepsilon_y\rightarrow\infty$, $||(\textbf{I}+\tau^\alpha\textbf{K})^{-1}||$ can be very close to zero;
(ii) if $\kappa_x$, $\kappa_y$ are larger than $\varepsilon_x$, $\varepsilon_y$ and $\varepsilon_x$, $\varepsilon_y$ themselves are small,
    $||(\textbf{I}+\tau^\alpha\textbf{K})^{-1}||$ can be larger than 1, however,  such issue can be remedied by increasing the grid numbers;
(iii) in general, the larger $M_x$, $M_y$, the smaller $||(\textbf{I}+\tau^\alpha\textbf{K})^{-1}||$;
(iv) when the computational domain expands,  $||(\textbf{I}+\tau^\alpha\textbf{K})^{-1}||$ grows at a speed, which may result in
    the invalidation of (\ref{ez83}) if $\varepsilon_x$, $\varepsilon_y$ and $M_x, M_y$ remain unchanged;
(v) when $\tau\rightarrow0$, the critical ratio between $\kappa_x, \kappa_y$ and
    $\varepsilon_x, \varepsilon_y$ to maintain this assumption appears to decrease, but it would be enhanced as the spatial grid is refined.

Consequently,  the assumption is meaningful and essentially a mild theoretical restriction in practise.

\section{Description of cubic B-spline DQ method}\label{s4}
In this section, a robust DQ method (MCB-DQM) based on the modified cubic B-splines $\{MB_m(x)\}_{m=0}^{M}$
is established for Eqs. (\ref{xz01})-(\ref{xz03}) by introducing the DQ approximations to fractional derivatives.
In the light of the essence of traditional DQ methods, we consider
\begin{align}
&\frac{\partial ^{\beta_1}u(x_i,y_j,t)}{\partial x^{\beta_1}}\cong\sum\limits_{m=0}^{M_x} {a_{im}^{(\beta_1)}u(x_m,y_j,t)}, \label{xz07}\\
&\frac{\partial ^{\beta_2}u(x_i,y_j,t)}{\partial y^{\beta_2}}\cong\sum\limits_{m=0}^{M_y} {b_{jm}^{(\beta_2)}u(x_i,y_m,t)},\label{xz08}
\end{align}
for fractional derivatives  in constructing our DQ algorithm, where $0\leq i\leq M_x$, $0\leq j\leq M_y$ and the weighted coefficients
$a_{im}^{(\beta_1)}$, $b_{jm}^{(\beta_2)}$ satisfy
\begin{small}
\begin{align}
\frac{\partial^{\beta_1} M{B_k}(x_i)}{\partial x^{\beta_1}}=\sum^{M_x}_{m=0}a_{im}^{(\beta_1)}M{B_k}(x_m),\ \ 0\leq i,k\leq M_x,\label{xz17}\\
\frac{\partial^{\beta_2} M{B_k}(y_j)}{\partial y^{\beta_2}}=\sum^{M_y}_{m=0}b_{jm}^{(\beta_2)}M{B_k}(y_m),\ \ 0\leq j,k\leq M_y.\label{xz18}
\end{align}
\end{small}
The validation of Eqs. (\ref{xz07})-(\ref{xz08}) is ensured by the linear properties of fractional derivatives.
$a_{im}^{(\beta_1)}$, $b_{jm}^{(\beta_2)}$ are then determined by tackling the resulting algebraic problems from
the above equations for each axis if the values of $\{MB_m(x)\}_{m=0}^{M}$ at all sampling points are known.

\subsection{The explicit formulas of fractional derivatives}
It is the weakly singular integral structure that makes it difficult to calculate the values of the fractional derivatives for a function as B-spline
at a sampling point. In the text that follows, we concentrate on the explicit expressions of the $\beta$-th ($1<\beta\leq2$) Riemann-Liouville
derivative of $\{B_m(x)\}_{m=-1}^{M+1}$ with a recursive technique. Since these basis splines are piecewise and locally compact on four consecutive
subintervals, we have
\begin{align*}
{^{RL}_{x_0}}D^\beta_xB_m(x)\!=\!\left\{ \begin{array}{l}
0,\qquad\qquad\qquad\ \,  x \in [x_0,x_{m-2})\\
{^{RL}_{x_{m-2}}}D^\beta_x\varphi_1(x), \qquad x \in [{x_{m - 2}},{x_{m - 1}})\\
{^{RL}_{x_{m-2}}}D^\beta_{x_{m-1}}\varphi_1(x)\\
 \quad +{^{RL}_{x_{m-1}}}D^\beta_x\varphi_2(x), \ \, x \in [{x_{m - 1}},{x_m})\\
{^{RL}_{x_{m-2}}}D^\beta_{x_{m-1}}\varphi_1(x)\\
 \quad +{^{RL}_{x_{m-1}}}D^\beta_{x_m}\varphi_2(x)\\
 \quad +{^{RL}_{x_m}}D^\beta_{x}\varphi_3(x), \quad \ x \in [{x_m},{x_{m + 1}})\\
{^{RL}_{x_{m-2}}}D^\beta_{x_{m-1}}\varphi_1(x)\\
 \quad +{^{RL}_{x_{m-1}}}D^\beta_{x_m}\varphi_2(x)\\
 \quad +{^{RL}_{x_m}}D^\beta_{x_{m+1}}\varphi_3(x)\\
 \quad +{^{RL}_{x_{m+1}}}D^\beta_{x}\varphi_4(x),  \ \ x \in [{x_{m + 1}},{x_{m + 2}})\\
{^{RL}_{x_{m-2}}}D^\beta_{x_{m-1}}\varphi_1(x)\\
 \quad +{^{RL}_{x_{m-1}}}D^\beta_{x_m}\varphi_2(x)\\
 \quad +{^{RL}_{x_m}}D^\beta_{x_{m+1}}\varphi_3(x)\\
 \quad +{^{RL}_{x_{m+1}}}D^\beta_{x_{m+2}}\varphi_4(x),  \ \, x \in [{x_{m + 2}},{x_M}]
\end{array} \right.
\end{align*}
with $2\leq m\leq M-2$. The compact supports of $B_{M-1}(x)$, $B_{M}(x)$, and $B_{M+1}(x)$ partially locate on the outside of $[x_0,x_M]$, so
do $B_{-1}(x)$, $B_{0}(x)$, and $B_{1}(x)$; nevertheless, $B_{M-1}(x)$, $B_{M}(x)$, and $B_{M+1}(x)$ can be thought of as
the special cases of the aforementioned argument, so are omitted here. Further, we have
\begin{align*}
{^{RL}_{x_0}}D^\beta_xB_{-1}(x)\!=\!\left\{ \begin{array}{l}
{^{RL}_{x_0}}D^\beta_x\varphi_4(x), \quad \  x \in [x_0,x_1)\\
{^{RL}_{x_0}}D^\beta_{x_1}\varphi_4(x), \quad \ x \in [x_1,{x_M}]
\end{array} \right.
\end{align*}
\begin{align*}
{^{RL}_{x_0}}D^\beta_xB_0(x)\!=\!\left\{ \begin{array}{l}
{^{RL}_{x_0}}D^\beta_x\varphi_3(x), \quad   x \in [x_0,x_1)\\
{^{RL}_{x_0}}D^\beta_{x_1}\varphi_3(x) \\
\quad +{^{RL}_{x_1}}D^\beta_x\varphi_4(x), \ \ \ x \in [x_1,x_2)\\
{^{RL}_{x_0}}D^\beta_{x_1}\varphi_3(x) \\
\quad +{^{RL}_{x_1}}D^\beta_{x_2}\varphi_4(x), \ \  x \in [x_2,{x_M}]
\end{array} \right.
\end{align*}
 \begin{align*}
{^{RL}_{x_0}}D^\beta_xB_1(x)\!=\!\left\{ \begin{array}{l}
{^{RL}_{x_0}}D^\beta_x\varphi_2(x), \quad   x \in [x_0,x_1)\\
{^{RL}_{x_0}}D^\beta_{x_1}\varphi_2(x) \\
\quad +{^{RL}_{x_1}}D^\beta_x\varphi_3(x), \ \ \ x \in [x_1,x_2)\\
{^{RL}_{x_0}}D^\beta_{x_1}\varphi_2(x) \\
\quad +{^{RL}_{x_1}}D^\beta_{x_2}\varphi_3(x)\\
\quad +{^{RL}_{x_2}}D^\beta_{x}\varphi_4(x), \ \  x \in [x_2,x_3)\\
{^{RL}_{x_0}}D^\beta_{x_1}\varphi_2(x) \\
\quad +{^{RL}_{x_1}}D^\beta_{x_2}\varphi_3(x)\\
\quad +{^{RL}_{x_2}}D^\beta_{x_3}\varphi_4(x). \ \  x \in [x_3,{x_M}]
\end{array} \right.
\end{align*}
On the other hand, as the integrands of the integration in fractional derivatives, $\varphi_i(x)$, $i=1,2,3,4$,
are cubic polynomials, for which, the order shrinks by one each time integration by parts is applied. Being aware of this,
we can eliminate the weakly singular integrations by repeating integration by parts four times for each $\varphi_i(x)$
to derive the fully explicit formulas. The derivation processes are lengthy and tedious, we therefore
outline the specific expressions of ${^{RL}_{x_0}}D^\beta_xB_m(x)$ in \textbf{Appendix}.

\subsection{Construction of cubic B-spline DQ method}
Use the early notations for brevity. On using DQ approximations (\ref{xz07})-(\ref{xz08}) to handle fractional derivatives,
Eq. (\ref{xz01}) is transformed into a set of first-order ODEs
\begin{align}
\begin{aligned}\label{xz09}
&\frac{\partial u(x_i,y_j,t)}{\partial t}-\varepsilon_x\sum\limits_{m=0}^{M_x}{a_{im}^{(\beta_1)}u(x_m,y_j,t)}\\
&\qquad \ -\varepsilon_y\sum\limits_{m=0}^{M_y}{b_{jm}^{(\beta_2)}u(x_i,y_m,t)}=f(x_i,y_j,t),
\end{aligned}
\end{align}
with $i=0,1,\cdots,M_x$, $j=0,1,\cdots,M_y$. Imposing the boundary constraint (\ref{xz03}) and
applying the Crank-Nicolson scheme, we thus obtain the following DQ scheme
\begin{align}
\left\{
\begin{aligned}\label{xz12}
&U_{ij}^n-\frac{\tau\varepsilon_x}{2}\sum_{m=1}^{M_x-1}a_{im}^{(\beta_1)}U_{mj}^n-\frac{\tau\varepsilon_y}{2}\sum_{m=1}^{M_y-1}b_{jm}^{(\beta_2)}U_{im}^n\\
& =U_{ij}^{n-1}+\frac{\tau\varepsilon_x}{2}\sum_{m=1}^{M_x-1}a_{im}^{(\beta_1)}U_{mj}^{n-1} \\
&\qquad\quad+\frac{\tau\varepsilon_y}{2}\sum_{m=1}^{M_y-1}b_{jm}^{(\beta_2)}U_{im}^{n-1}+\tau f^{n-1/2}_{ij},
\end{aligned}\right.
\end{align}
where $i=1,2,\cdots,M_x-1$, $j=1,2,\cdots,M_y-1$. It is visible that DQ methods are truly meshless and convenient in implementation.
Due to the insensitivity to dimensional changes, (\ref{xz12}) can easily be generalized to the higher-dimensional space-fractional problems,
but do not cause the rapid increase of computing burden.

\section{Illustrative examples}\label{s5}
In this section, a couple of numerical examples are carried out to
gauge the practical performance of MCTB-DQM and new MCB-DQM. In order to check their accuracy,
we compute the errors by using the norms
\begin{align*}
&e_\infty(M)\cong\max_{i}\Big|u^n_i-U^n_i\Big|,\\
&e_2(M)\cong\sqrt{\frac{1}{M}\sum^{M-1}_{i=1}\Big|u^n_i-U^n_i\Big|^2},\\
&e_N(M)\cong\sqrt{\sum^{M-1}_{i=1}\Big|u^n_{i}-U^n_{i}\Big|^2\bigg/\sum^{M-1}_{i=1}\Big|U_{i}^0\Big|^2},\\
&e_\infty(M_x,M_y)\cong\max_{i,j}\Big|u^n_{ij}-U^n_{ij}\Big|,\\
&e_2(M_x,M_y)\cong\sqrt{\frac{1}{M_xM_y}\sum^{M_x-1}_{i=1}\sum^{M_y-1}_{j=1}\Big|u^n_{ij}-U^n_{ij}\Big|^2},
\end{align*}
where $e_N(M)$ is termed by a normalized $L_2$-norm.
As to $\{\omega^\alpha_k\}_{k=0}^{n}$ in Eqs. (\ref{eq22})-(\ref{eq23}), we use (\ref{ez06}) in the first and fifth examples
and (\ref{eq08}) in the others but not the last two ones.
In the computation, our algorithms are implemented on Matlab platform in a Lenovo PC with Intel(R) Pentium(R) G2030 3.00GHz CPU and 4 GB RAM.
The obtained results are comparatively discussed with the early works available in the open literature.\\

\noindent
\emph{Example 6.1.} Let $\kappa=1$, $\varepsilon=2$; the Eqs. (\ref{eq01})-(\ref{eq03})
with $\psi(x)=\exp(x)$, $g_1(t)=E(t^\alpha)$, $g_2(t)=eE(t^\alpha)$
and homogeneous forcing term are considered on $[0,1]$, where $E(t^\alpha)$ is the well-known
\emph{Mittag-Leffler function}
\begin{equation*}
    E_\alpha(z)=\sum_{k=0}^{\infty}\frac{z^k}{\Gamma(\alpha k+1)}, \quad 0<\alpha<1.
\end{equation*}
It is verified that its solution is $u(x,t)=\exp(x)E(t^\alpha)$. In order to show the
convergence of MCTB-DQM, we fix $\tau=1.0\times10^{-5}$ so that the temporal errors are negligible as compared
to spatial errors. The numerical results at $t=0.1$ for various $\alpha$ are displayed
in Table \ref{tab1}; the convergent rate is shortly written as ``Cov. rate''. As one sees, our method
is pretty stable and convergent with almost spatial second-order for this problem.\\

\noindent
\emph{Example 6.2.} In this test, we solve a diffusion equation on $[0,1]$ with $\varepsilon=1$,
$\psi(x)=4x(1-x)$, zero boundary condition and right side. Its true solution has the form 
\begin{equation*}
  u(x,t)=\frac{16}{\pi^3}\sum^\infty_{k=1}\frac{1}{k^3}E_\alpha(-k^2\pi^2t^\alpha)(1-(-1)^k)\sin(k\pi x).
\end{equation*}
For comparison of the numerical results given by FDS-D I, FDS-D II \cite{28} and the semi-discrete FEM \cite{29},
we choose the same time stepsize $\tau=1.0\times10^{-4}$.
Letting $\alpha=0.1$, $0.5$ and $0.95$, the corresponding results of these four methods at $t=1$ are tabulated
side by side in Table \ref{tab2}, from which, we conclude that  MCTB-DQM is accurate and produces very small errors as
the other three methods as the grid number $M$ increases. \\  

\begin{table*}
\centering
\caption{The numerical results at $t=0.1$ with $\tau=1.0\times10^{-5}$ for Example 6.1.} \label{tab1}
\begin{tabular}{lllclc}
\toprule
          $\alpha$  &  $M$   &$e_2(M)$ &Cov. rate &$e_\infty(M)$ & Cov. rate  \\
\midrule $0.2$      & 8        &2.4430e-03  &-       &    3.5200e-03  &-     \\
                    & 16       &6.3696e-04  &1.9394  &    9.2142e-04  &1.9337   \\
                    & 32       &1.6272e-04  &1.9688  &    2.4362e-04  &1.9192   \\
                    & 64       &4.1425e-05  &1.9738  &    6.2649e-05  &1.9592   \\
                    & 128      &1.0765e-05  &1.9441  &    1.5906e-05  &1.9777   \\
       $0.5$        & 8        &1.2489e-03  &-       &    1.8283e-03  &-     \\
                    & 16       &3.2655e-04  &1.9352  &    4.8198e-04  &1.9235   \\
                    & 32       &8.3679e-05  &1.9644  &    1.2771e-04  &1.9160   \\
                    & 64       &2.1466e-05  &1.9628  &    3.3008e-05  &1.9520   \\
                    & 128      &5.7295e-06  &1.9056  &    8.4114e-06  &1.9724   \\
       $0.8$        & 8        &9.7261e-04  &-       &    1.4452e-03  &-     \\
                    & 16       &2.5487e-04  &1.9321  &    3.8329e-04  &1.9147   \\
                    & 32       &6.5378e-05  &1.9629  &    1.0160e-04  &1.9156   \\
                    & 64       &1.6774e-05  &1.9625  &    2.6330e-05  &1.9481  \\
                    & 128      &4.4723e-06  &1.9072  &    6.7183e-06  &1.9705   \\
\bottomrule
\end{tabular}
\end{table*}

\begin{table*}
\centering
\caption{A comparison of $e_N(M)$ at $t=1$ with $\tau=1.0\times10^{-4}$ for Example 6.2} \label{tab2}
\begin{tabular}{clcccc}
\toprule
          $\alpha$  &  $M$   &FDS-D I \cite{28} &FDS-D II \cite{28} &FEM \cite{29} &MCTB-DQM  \\
\midrule $0.1$      & 8        &9.98e-04  &1.00e-03  &    5.23e-04  &2.7445e-04     \\
                    & 16       &2.44e-04  &2.53e-04  &    1.29e-04  &3.6823e-05   \\
                    & 32       &5.36e-05  &6.33e-05  &    3.21e-05  &5.0524e-06   \\
                    & 64       &5.89e-06  &1.55e-05  &    8.01e-06  &9.6533e-07   \\
                    & 128      &6.08e-06  &3.62e-06  &    2.00e-06  &4.6423e-07   \\
       $0.5$        & 8        &7.13e-04  &7.13e-04  &    3.37e-04  &1.9258e-04     \\
                    & 16       &1.79e-04  &1.79e-04  &    8.31e-04  &2.6563e-05   \\
                    & 32       &4.46e-05  &4.44e-05  &    2.07e-05  &4.5588e-06   \\
                    & 64       &1.07e-05  &1.06e-05  &    5.17e-06  &1.7866e-06   \\
                    & 128      &2.23e-06  &2.12e-06  &    1.30e-06  &1.4534e-06   \\
       $0.95$        & 8       &1.11e-04  &1.11e-04  &    4.84e-05  &2.8976e-05     \\
                    & 16       &2.83e-05  &2.82e-05  &    1.21e-05  &4.0937e-06   \\
                    & 32       &7.09e-06  &7.08e-06  &    3.05e-06  &8.7766e-07   \\
                    & 64       &1.76e-06  &1.75e-06  &    7.93e-07  &4.7682e-07  \\
                    & 128      &4.29e-07  &4.23e-07  &    2.32e-07  &4.2766e-07   \\
\bottomrule
\end{tabular}
\end{table*}

\begin{figure*}
\centering
\subfigure{\includegraphics[width=2.85in]{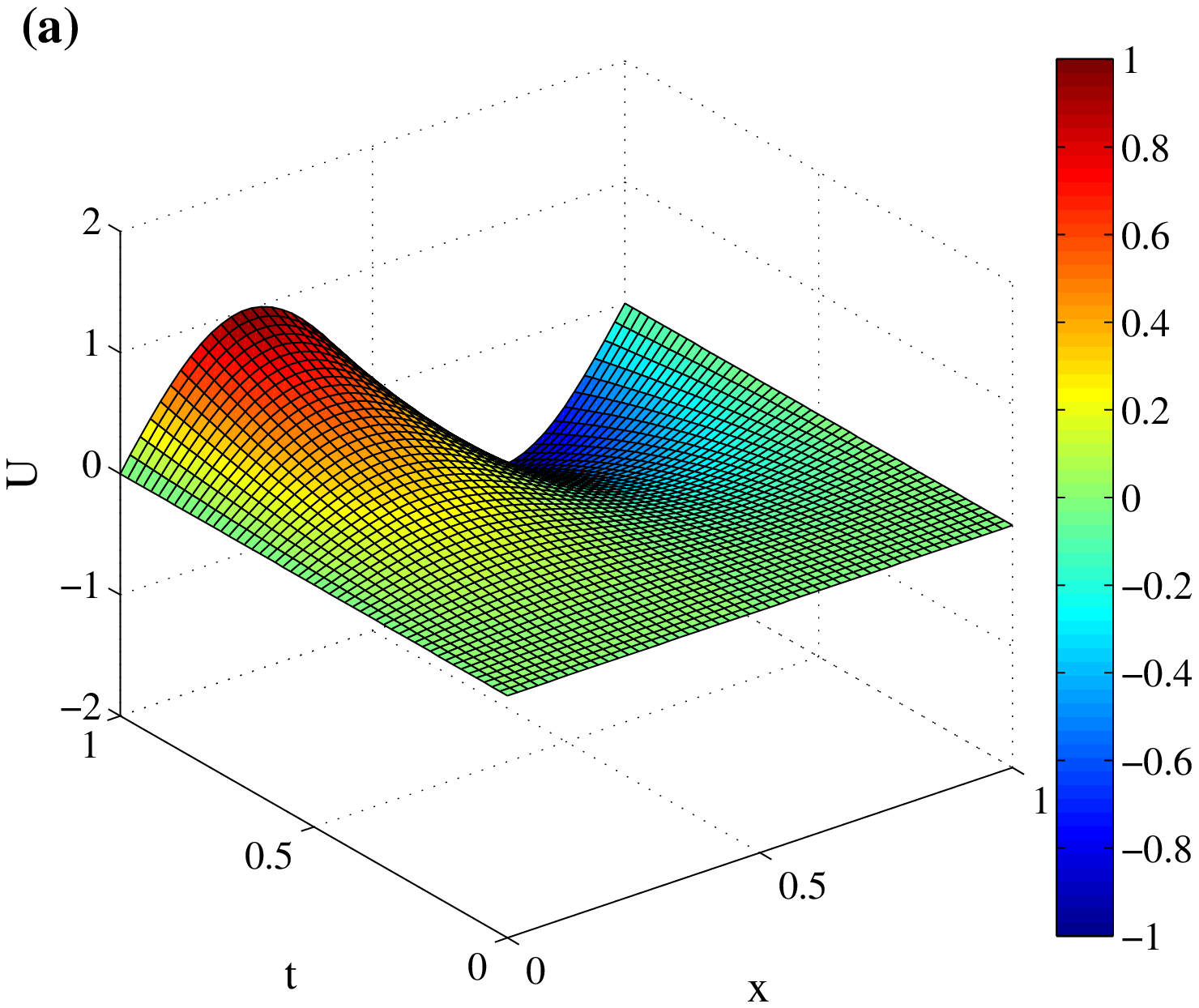}}\hspace{0.28in}
\subfigure{\includegraphics[width=2.85in]{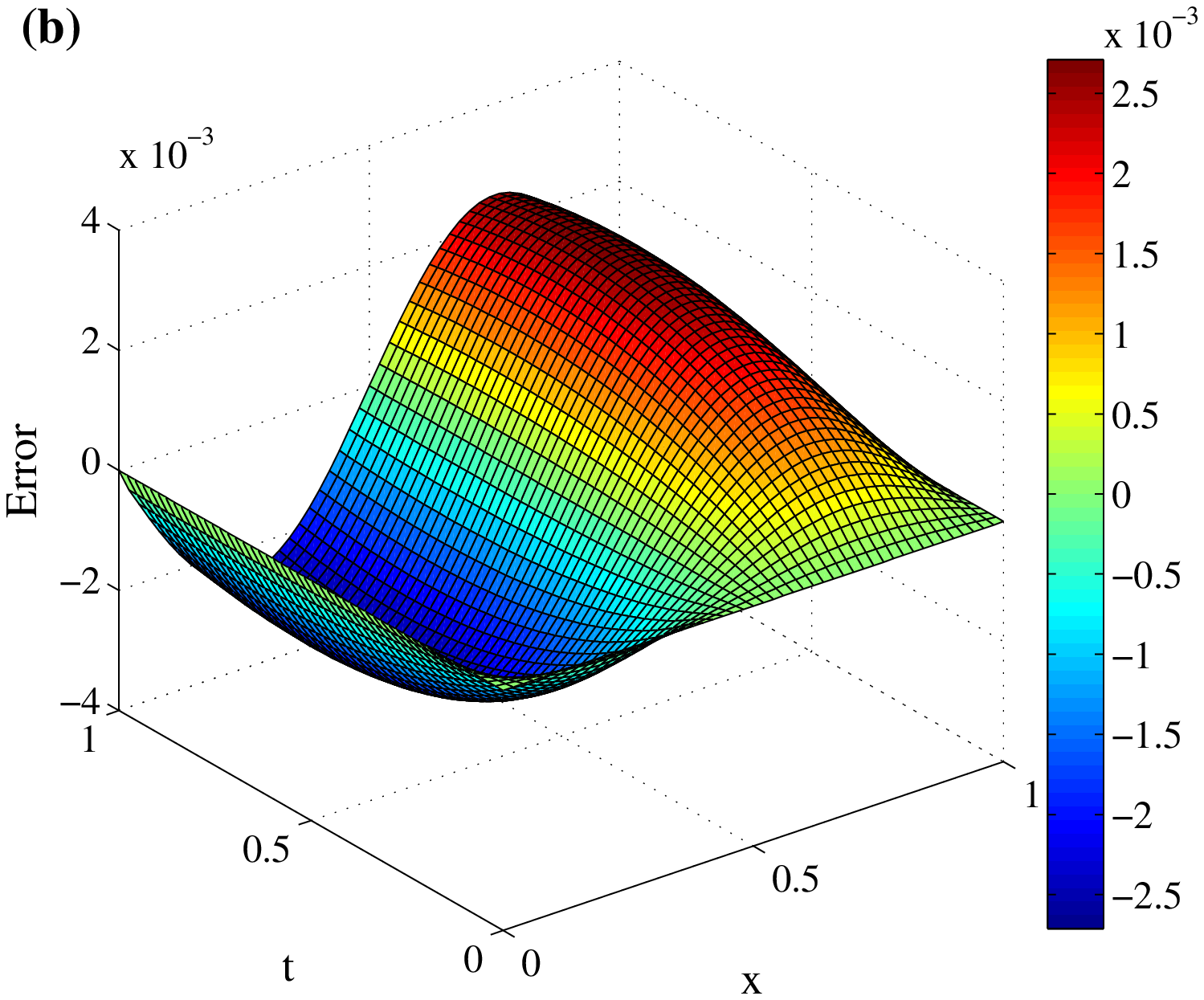}}
\caption{The approximate solution and error distribution at $t=1$ with $\alpha=0.8$ for Example 6.3}\label{fig1}
\end{figure*}

\noindent
\emph{Example 6.3.} Let $\kappa=0$, $\varepsilon=1$; we solve Eqs. (\ref{eq01})-(\ref{eq03})
with homogeneous initial and boundary values, and 
\begin{equation*}
   f(x,t)=\frac{2t^{2-\alpha}\sin(2\pi x)}{\Gamma(3-\alpha)}+4\pi^2t^2\sin(2\pi x),
\end{equation*}
on $[0,1]$. The ture solution is $u(x,t)=t^2\sin(2\pi x)$. The algorithm is first run with
$\alpha=0.8$, $\tau=2.0\times10^{-2}$ and $M=50$. In Fig. \ref{fig1}, we plot the approximate solution and a point to point
error distribution at $t=1$, where good accuracy is observed. In Table \ref{tab3}, we then report a comparison of $e_2(M)$,
$e_\infty(M)$ at $t=1$ between MCTB-DQM and CBCM \cite{ref23}, when $\alpha=0.3$. Here, MCTB-DQM uses $\tau=5.0\times10^{-3}$ while
CBCM chooses $\tau=1.25\times10^{-3}$. As expected, our approach generates the  approximate solutions
with a better accuracy than those obtained by CBCM. \\

\begin{table*}
\centering
\caption{A comparison of $e_2(M)$, $e_\infty(M)$  at $t=1$ with $\alpha=0.3$ for Example 6.3}\label{tab3}
\begin{tabular}{lllll}
\toprule
\multicolumn{1}{l}{\multirow{2}{0.6cm}{$M$}}
&\multicolumn{2}{l}{CBCM \cite{ref23}} &\multicolumn{2}{l}{MCTB-DQM} \\
\cline{2-5}& $e_2(M)$  &$e_\infty(M)$ &$e_2(M)$  &$e_\infty(M)$  \\
\midrule 8     &3.4134e-02  &4.8273e-02  &  9.4300e-03  &1.5762e-02  \\
        16     &8.7334e-03  &1.2351e-02  &  1.1924e-03  &2.1670e-03  \\
        32     &2.1955e-03  &3.1048e-03  &  1.5040e-04  &2.8541e-04  \\
        64     &5.4957e-04  &7.7721e-04  &  1.8925e-05  &3.6701e-05  \\
        128     &1.3739e-04  &1.9430e-04  &  2.3752e-06  &4.6559e-06  \\
\bottomrule
\end{tabular}
\end{table*}

\begin{figure*}
\centering
\subfigure{\includegraphics[width=2.7in]{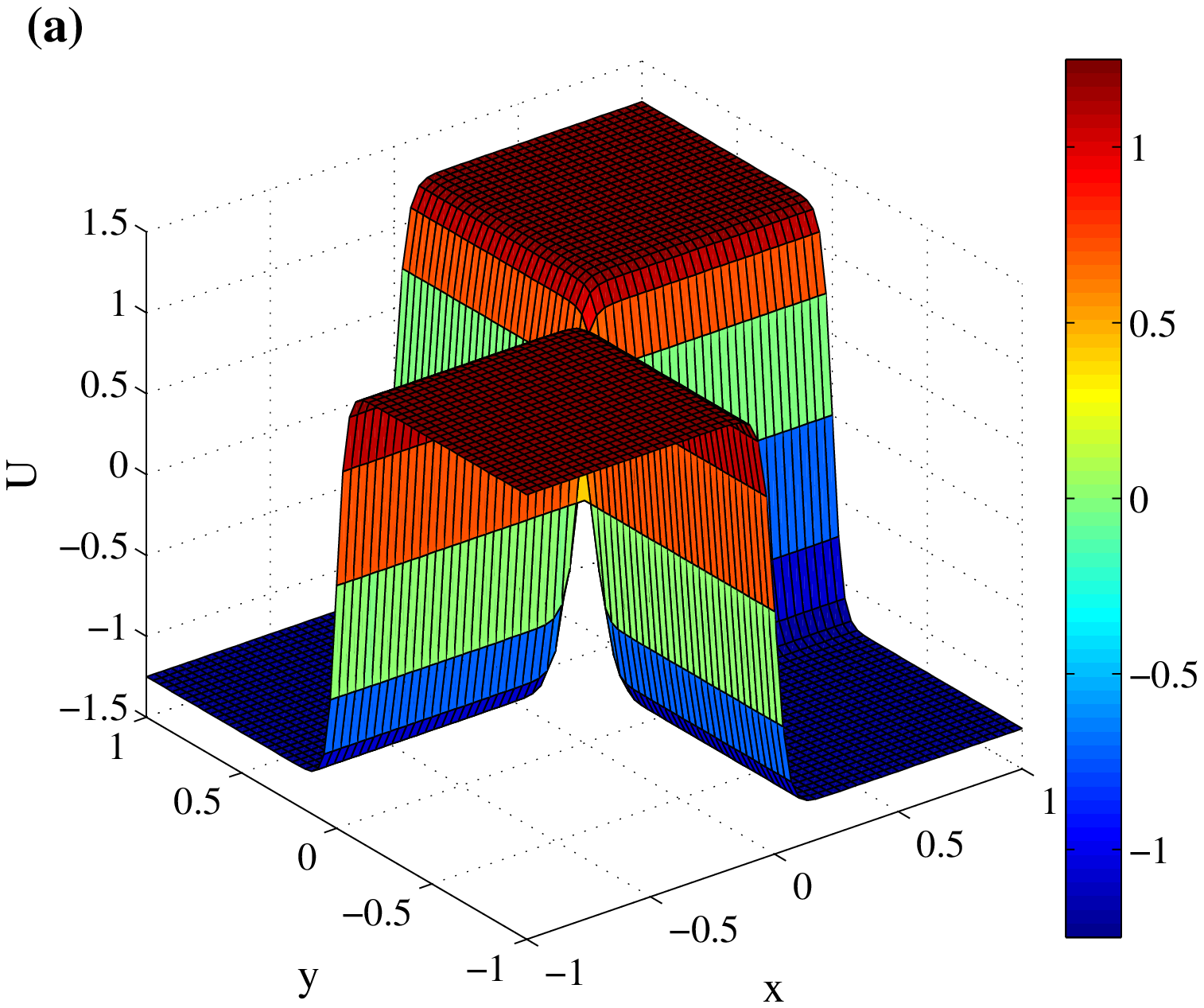}}\hspace{0.28in}
\subfigure{\includegraphics[width=2.7in]{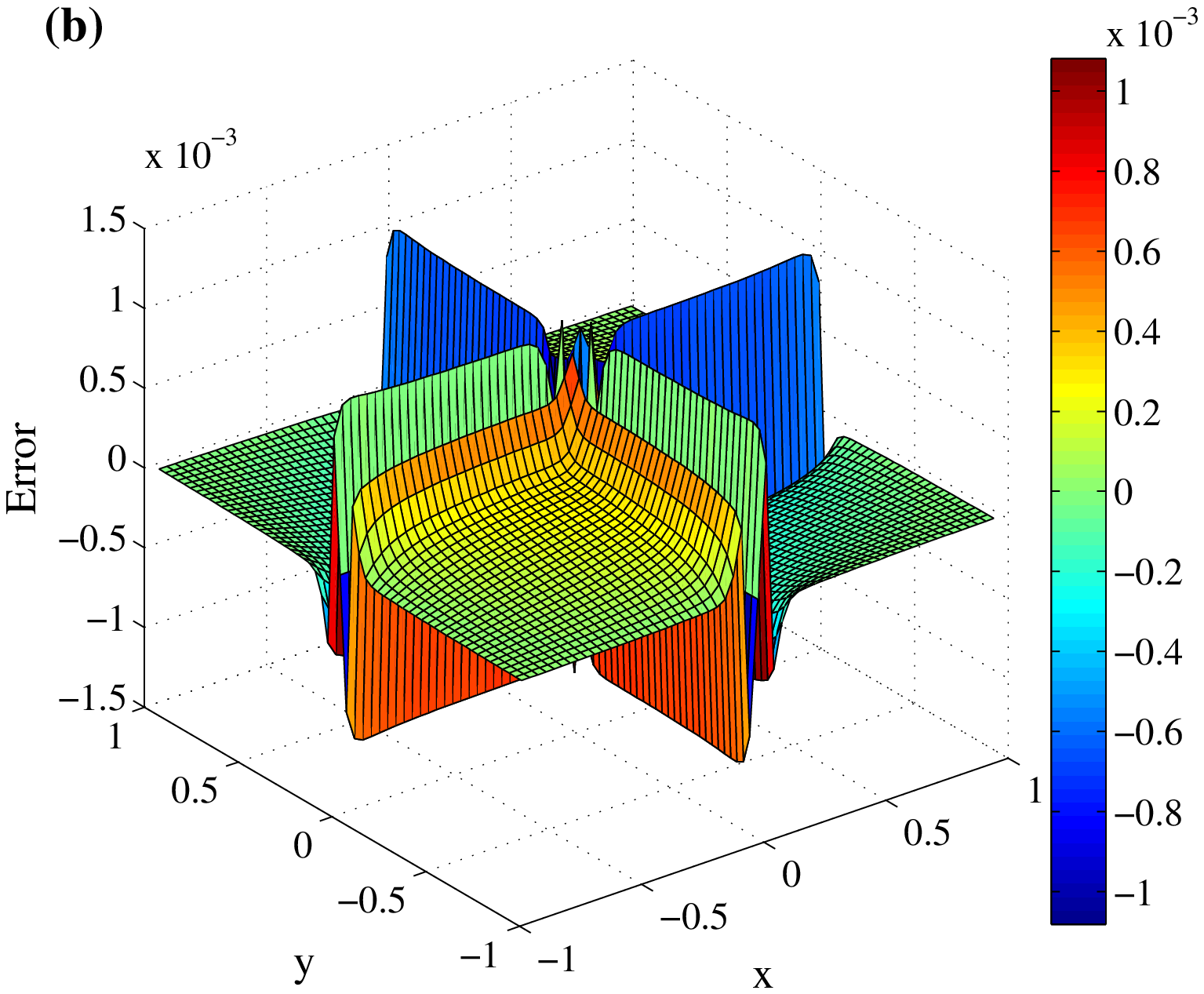}}
\caption{The approximate solution and error distribution at $t=0.5$ with $\alpha=0.5$ for Example 6.4}\label{fig2}
\end{figure*}

\noindent
\emph{Example 6.4.} We consider a 2D diffusion equation on $[-1,1]\times[-1,1]$ with
$\varepsilon_x=\varepsilon_y=1$, which is referred to by Zhai and Feng as a test of a block-centered
finite difference method (BCFDM) on nonuniform grids \cite{31}. The forcing function is specified to
enforce 
\begin{equation*}
    u(x,y,t)=(1+t^2)\tanh(20x)\tanh(20y).
\end{equation*}
Under $\tau=1.0\times10^{-2}$, $M_x=M_y=60$ and $\alpha=0.5$,  we first plot the approximate solution and
a point to point error distribution at $t=0.5$ in Fig. \ref{fig2}. Then, 
we compare MCTB-DQM and BCFDM in term of $e_\infty(M_x,M_y)$
at $t=0.5$ in Table \ref{tab4}. 
It is obvious that MCTB-DQM produces significantly smaller errors than BCFDM as the grid number increases
despite a smaller time stepsize $\tau=2.5\times10^{-3}$ and the nonuniform girds BCFDM adopts; moreover, MCTB-DQM
provides more than quadratic rate of convergence for this problem. \\

\begin{table*}
\centering
\caption{A comparison of $e_\infty(M_x,M_y)$ at $t=0.5$ with $\alpha=0.5$ for Example 6.4} \label{tab4}
\begin{tabular}{ccccc}
\toprule
\multicolumn{1}{c}{\multirow{2}{1.1cm}{$M_x, M_y$}}
&\multicolumn{2}{l}{BCFDM \cite{31}} &\multicolumn{2}{l}{MCTB-DQM} \\
\cline{2-5}&$e_\infty(M_x,M_y)$ &Cov. rate  &$e_\infty(M_x,M_y)$ & Cov. rate \\
\midrule12     &8.75e-02\ \  &-     &   3.3376e-01  &- \\
        24     &2.73e-02\ \  &1.68  &   4.6331e-03  &6.1707   \\
        48     &8.26e-03\ \  &1.73  &   3.4566e-04  &3.7446   \\
        96     &2.24e-03\ \  &1.88  &   1.8605e-05  &4.2156   \\
\bottomrule
\end{tabular}
\end{table*}

\begin{table*}[!htb]
\centering
\caption{The numerical results in term of $e_2(M)$ at $t=0.1$ for Example 6.5}\label{tab6} 
\begin{tabular}{ccccc}
\toprule
$\alpha$   &  Method &Real part  & Imaginary part & CPU time \\
\midrule 0.2     &scheme (\ref{eq22})  &7.4588e-04  &7.3396e-04 & 10.355 (s) \\
         0.5     &scheme (\ref{eq22})  &3.9328e-03  &3.7620e-03 & 10.819 (s) \\
         0.8     &scheme (\ref{eq22})  &4.5245e-03  &4.6714e-03 & 10.775 (s) \\
         1.0     &scheme (\ref{eq22})  &2.2229e-03  &2.2153e-03 & 1.7436 (s)  \\
                 &scheme (\ref{ez09})  &8.0954e-04  &8.1843e-04 & 0.1521 (s) \\
\bottomrule
\end{tabular}
\end{table*}

\begin{figure*}[!htb]
\centering
\subfigure{\includegraphics[width=2.9in]{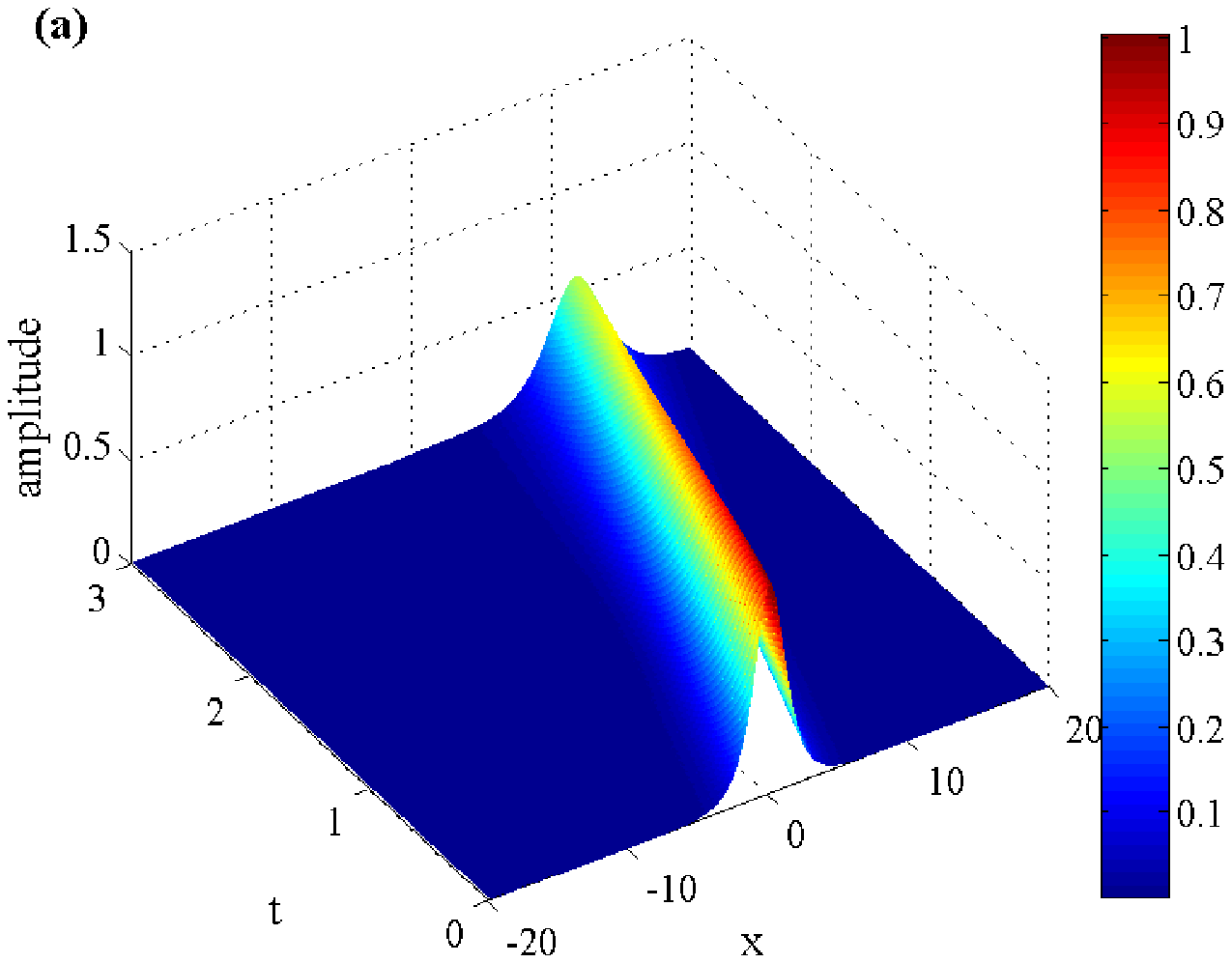}}\hspace{0.28in}
\subfigure{\includegraphics[width=2.9in]{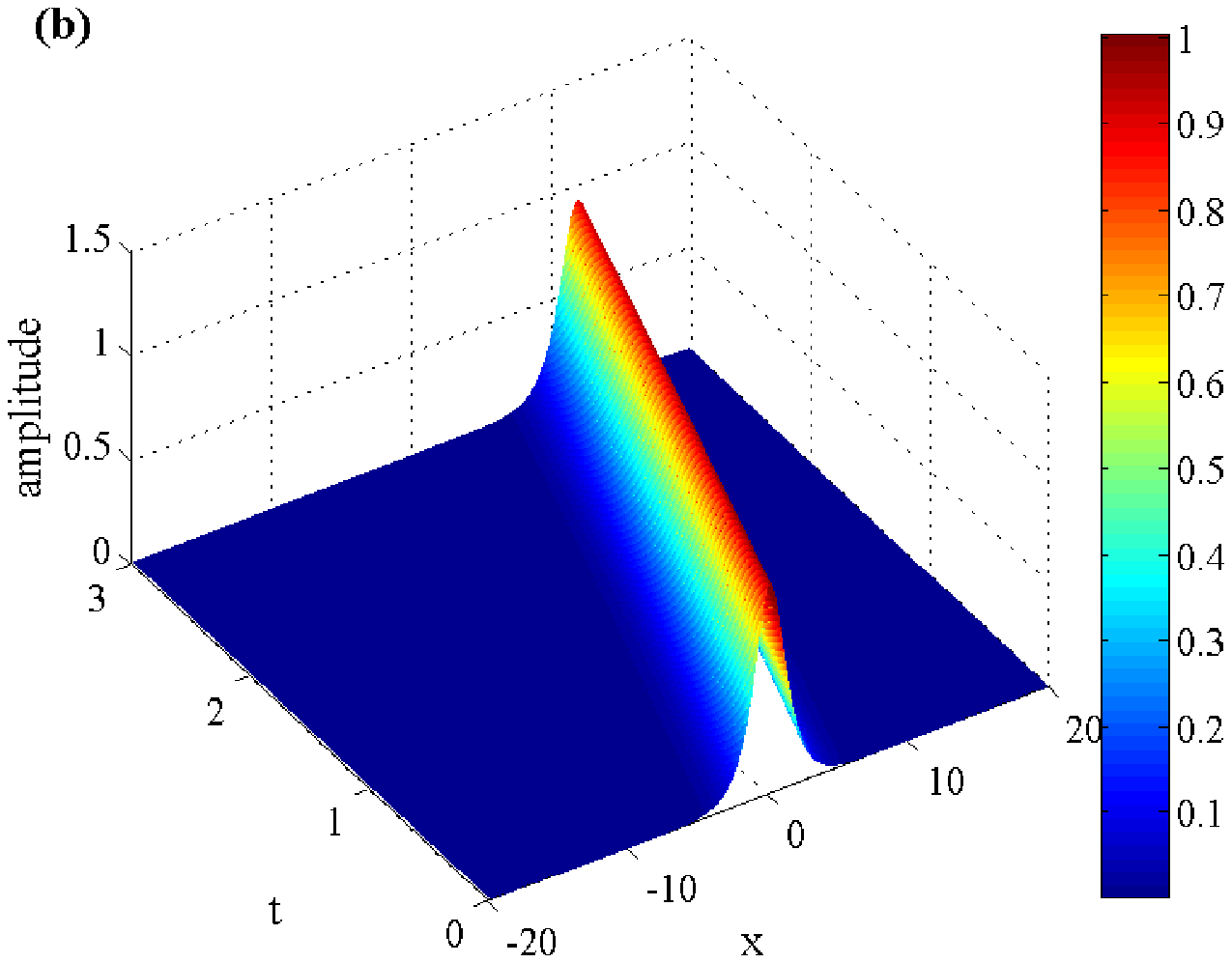}}
\caption{The single soliton propagation for $\alpha=0.98$, $1.0$ with $\tau=2.0\times10^{-3}$ and $M=200$.}\label{fig8}
\end{figure*}

\begin{figure*}[!htb]
\centering
\subfigure{\includegraphics[width=2.9in]{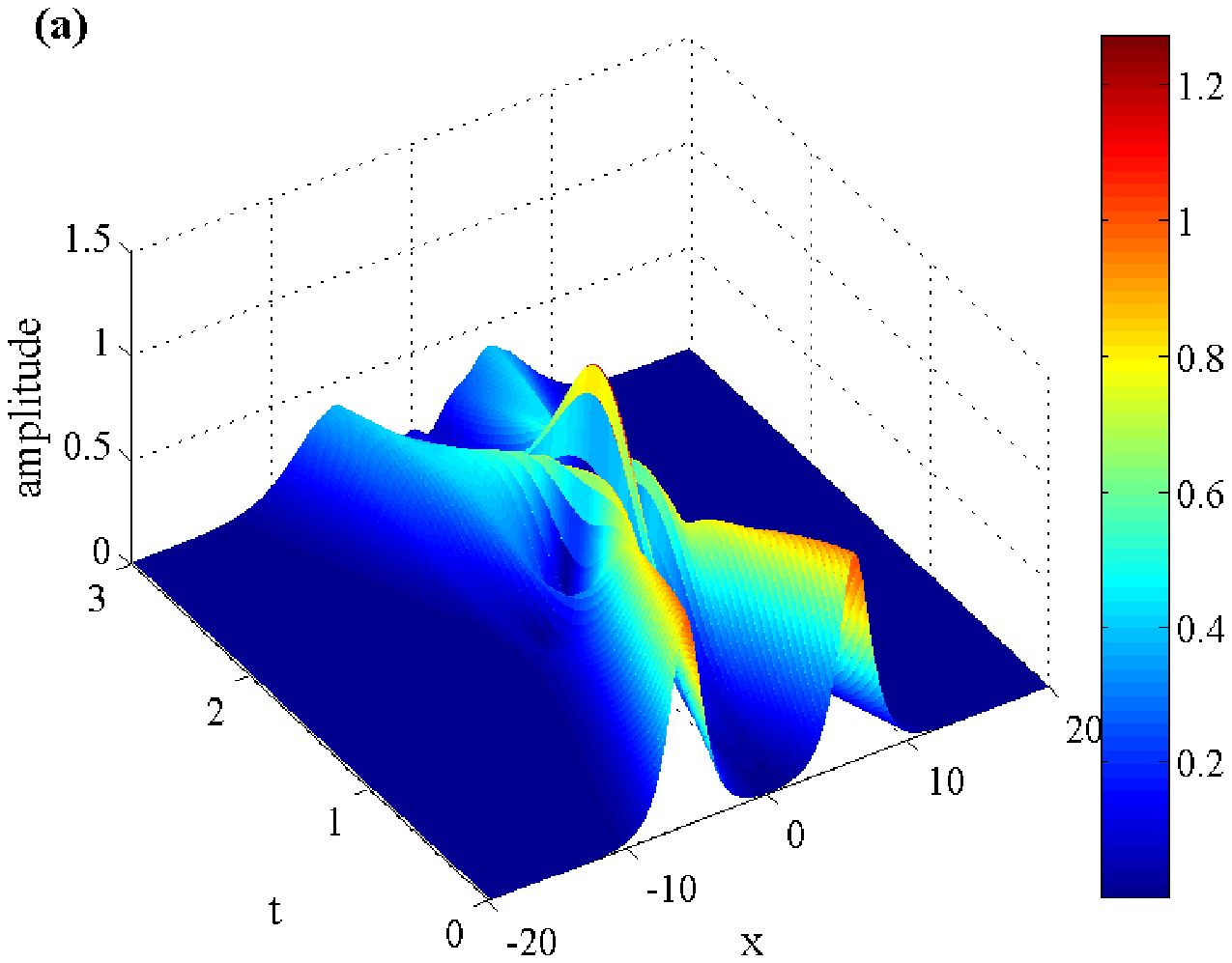}}\hspace{0.28in}
\subfigure{\includegraphics[width=2.9in]{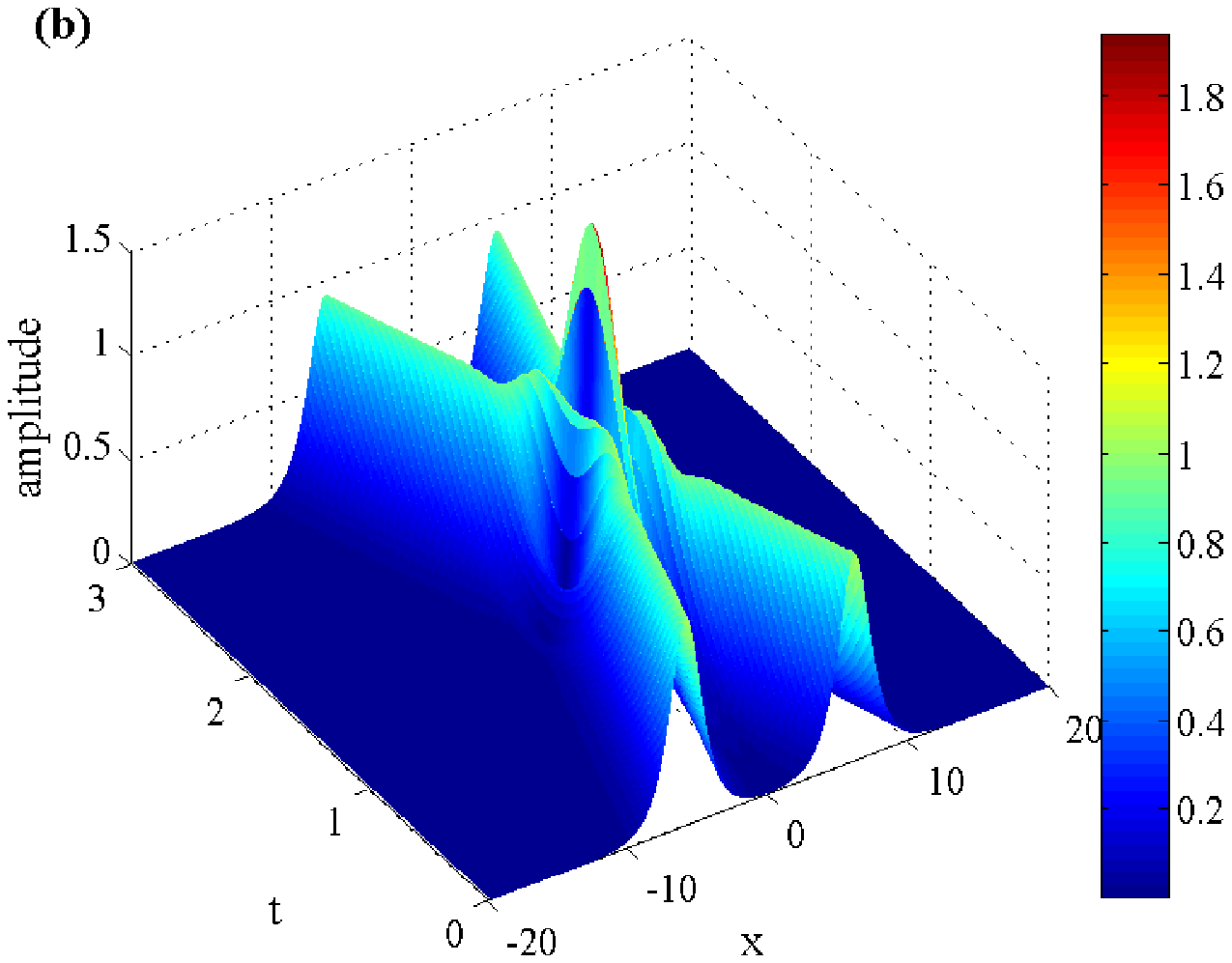}}
\caption{The interaction of double solitons for $\alpha=0.96$, $1.0$ with $\tau=2.0\times10^{-3}$ and $M=200$.}\label{fig9}
\end{figure*}

\noindent
\emph{Example 6.5.} In this test, we simulate the solitons propagation and collision governed by
    the following time-fractional nonlinear Schr\"{o}dinger equation (NLS):
\begin{align*}
    \textrm{i}\frac{\partial^\alpha u}{\partial t^\alpha}+\frac{\partial^2 u}{\partial x^2} + \beta|u|^2u=0, \ \ x\in(-\infty,+\infty),
\end{align*}
with $\textrm{i}=\sqrt{-1}$ and $\beta$ being a real constant, subjected to the initial values of two Gaussian types:
\begin{itemize}
  \item [(i)] mobile soliton
          \begin{equation}\label{ez90}
            \psi(x)=\textrm{sech}(x)\exp(2\textrm{i}x);
        \end{equation}
  \item [(ii)] double solitons collision
       \begin{equation}\label{ez91}
            \psi(x)=\sum_{j=1}^2\textrm{sech}(x-x_j)\exp(\textrm{i}p(x-x_j)).
       \end{equation}
\end{itemize}
When $\alpha=1$ and $\beta=2$, the NLS with Eq. (\ref{ez90}) has the soliton solution $u(x,t)=\textrm{sech}(x-4t)\exp(\textrm{i}(2x-3t))$.
As the solutions would generally decay to zero as $|x|\rightarrow \infty$,
we truncate the system into a bounded interval $\Omega=[a,b]$ with $a\ll0$ and $b\gg0$, and enforce periodic or homogeneous Dirichlet boundary conditions.
Letting $u(x,t)=U(x,t)+\textrm{i}V(x,t)$. Then, the original equation can be recast as a coupled diffusion system
\begin{align*}
  &\frac{\partial^\alpha U}{\partial t^\alpha}+\frac{\partial^2 V}{\partial x^2}+\beta(U^2+V^2)V=0, \\
  &\frac{\partial^\alpha V}{\partial t^\alpha}-\frac{\partial^2 U}{\partial x^2}-\beta(U^2+V^2)U=0.
\end{align*}
After applying the scheme (\ref{eq22}), nevertheless, a nonlinear system has to be solved at each time step. In such a case, the
Newton's iteration is utilized to treat it and terminated by reaching a solution with tolerant error $1.0\times10^{-12}$ if $\alpha=1$, for which,
the Jacobian matrix is
\begin{align*}
\textbf{J}=\left( \begin{array}{cc}
2UV&U^2+3V^2\\
-3U^2-V^2&-2UV
\end{array} \right).
\end{align*}
When $\alpha\neq1$, because the analytic solutions still remain unknown and the Newton's procedure relies heavily on its initial values,
we instead employ the trust-region-dogleg algorithm built into Matlab to improve the convergence of iteration.  At first, taking 
$\tau=2.0\times10^{-3}$, $M=100$, $\beta=2$, and $\Omega=[-10,10]$, the mean square errors at $t=0.1$ with the initial condition (\ref{ez90})
for various $\alpha$ are reported in Table \ref{tab6}, where the solutions computed by using the coefficients (\ref{eq08})
on a very fine time-space lattice, i.e., $\tau=2.5\times10^{-4}$, $M=400$,  are adopted as reference solutions ($\alpha\neq1$).
As seen from Table \ref{tab6}, our methods are convergent and applicable to nonlinear coupled problems; besides,
the scheme (\ref{ez09}) is clearly more efficient than (\ref{eq22}) since an extra Newton's outer loop is avoided.
Then, retaking $M=200$ and $\Omega=[-20,20]$, we display the evolution of the amplitude of the mobile soliton created  
by (\ref{eq22}) for $\alpha=0.98$ and $1.0$ in Fig. \ref{fig8}, respectively. Using the same discrete parameters, we consider
the double solitons collision for $\alpha=0.96$ and $1.0$ with $x_1=-6$, $x_2=6$, and $p=\pm2$ in Fig. \ref{fig9}.
It is easily drawn from these figures that the width and height of the solitons have been significantly changed by the fractional derivative.
In particular, when $\alpha=1$, a collision of double solitons without any reflection, transmission, trapping and creation of new solitary waves is exhibited,
which says that it is elastic, while in fractional cases, the shapes of the solitons may not be retained after they intersect each other. \\

\begin{figure*}
\centering
\subfigure{\includegraphics[width=2.7in]{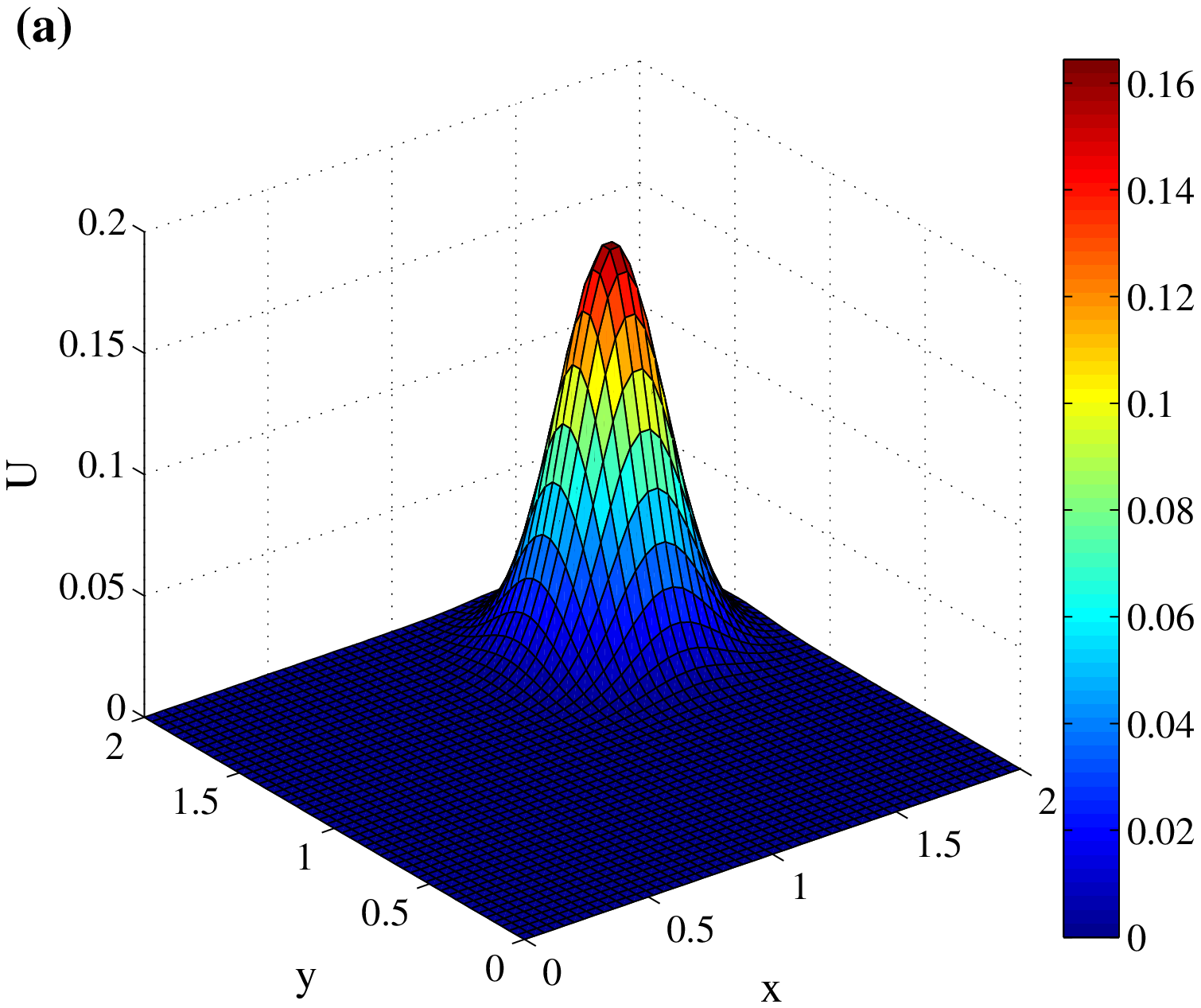}}\hspace{0.28in}
\subfigure{\includegraphics[width=2.7in]{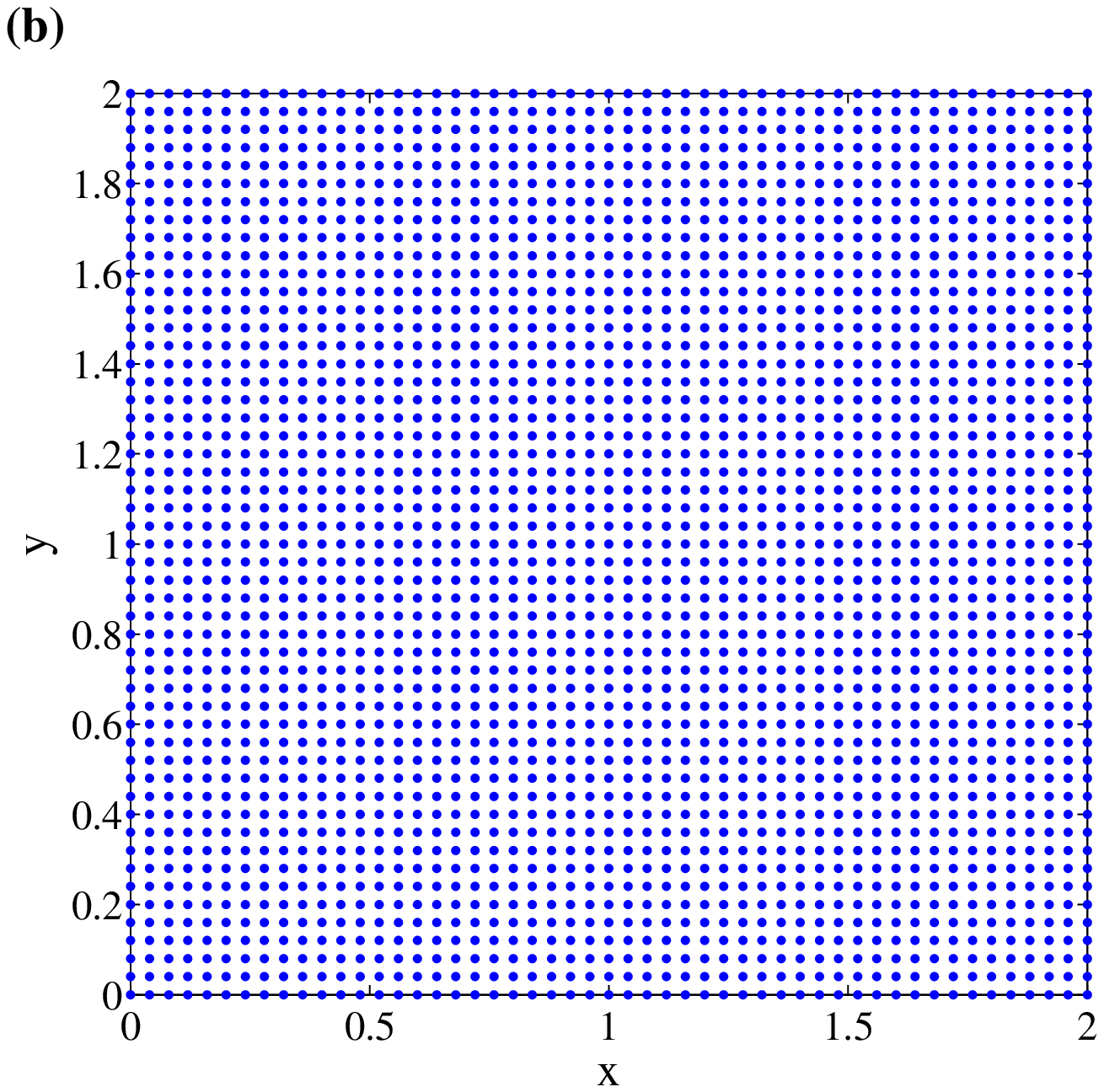}}
\caption{The true solution at $t=1.25$ and spatial lattice points for Example 6.6}\label{fig3}
\end{figure*}

\begin{figure*}
\centering
\subfigure{\includegraphics[width=2.7in]{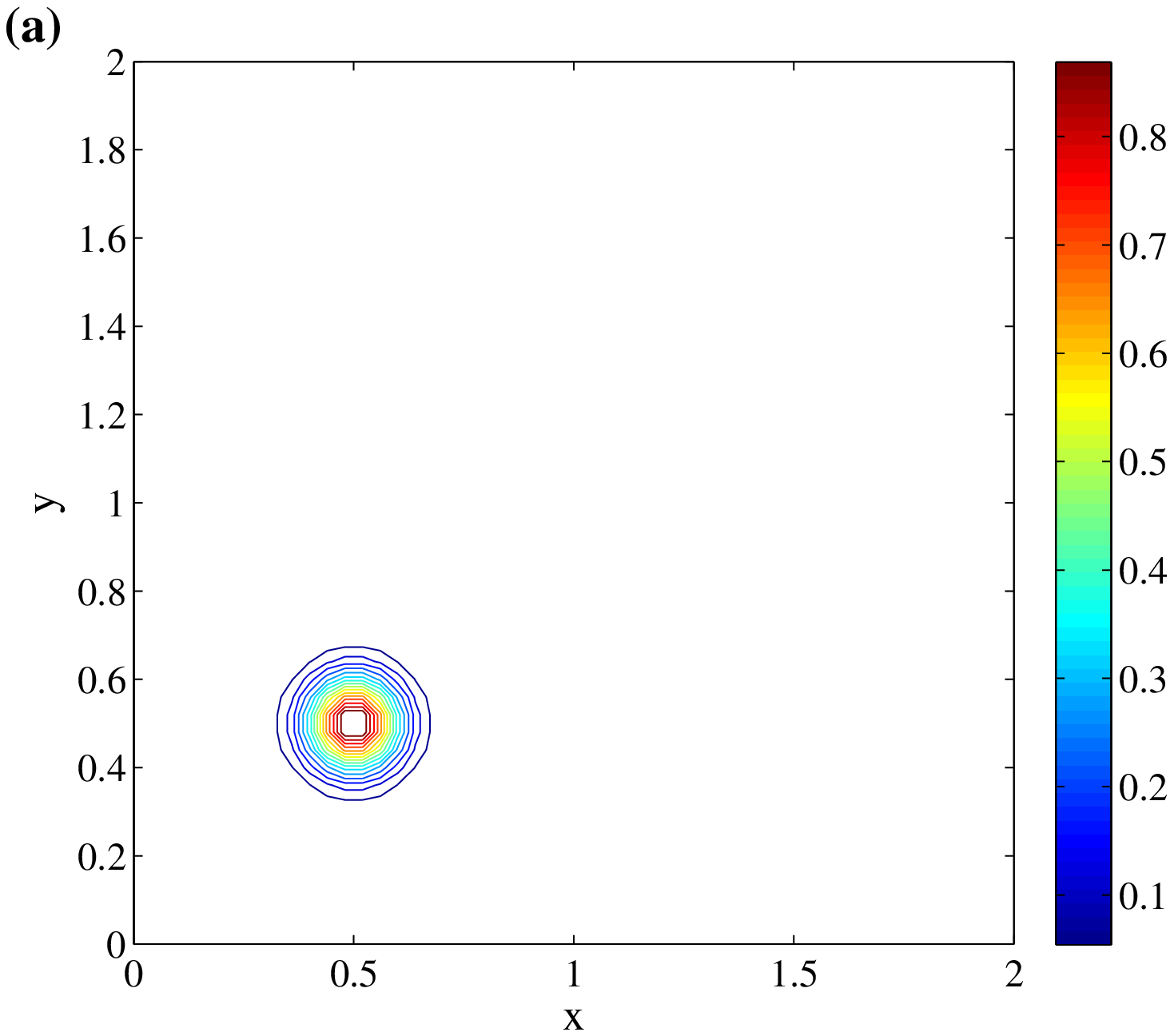}}\hspace{0.2in}
\subfigure{\includegraphics[width=2.7in]{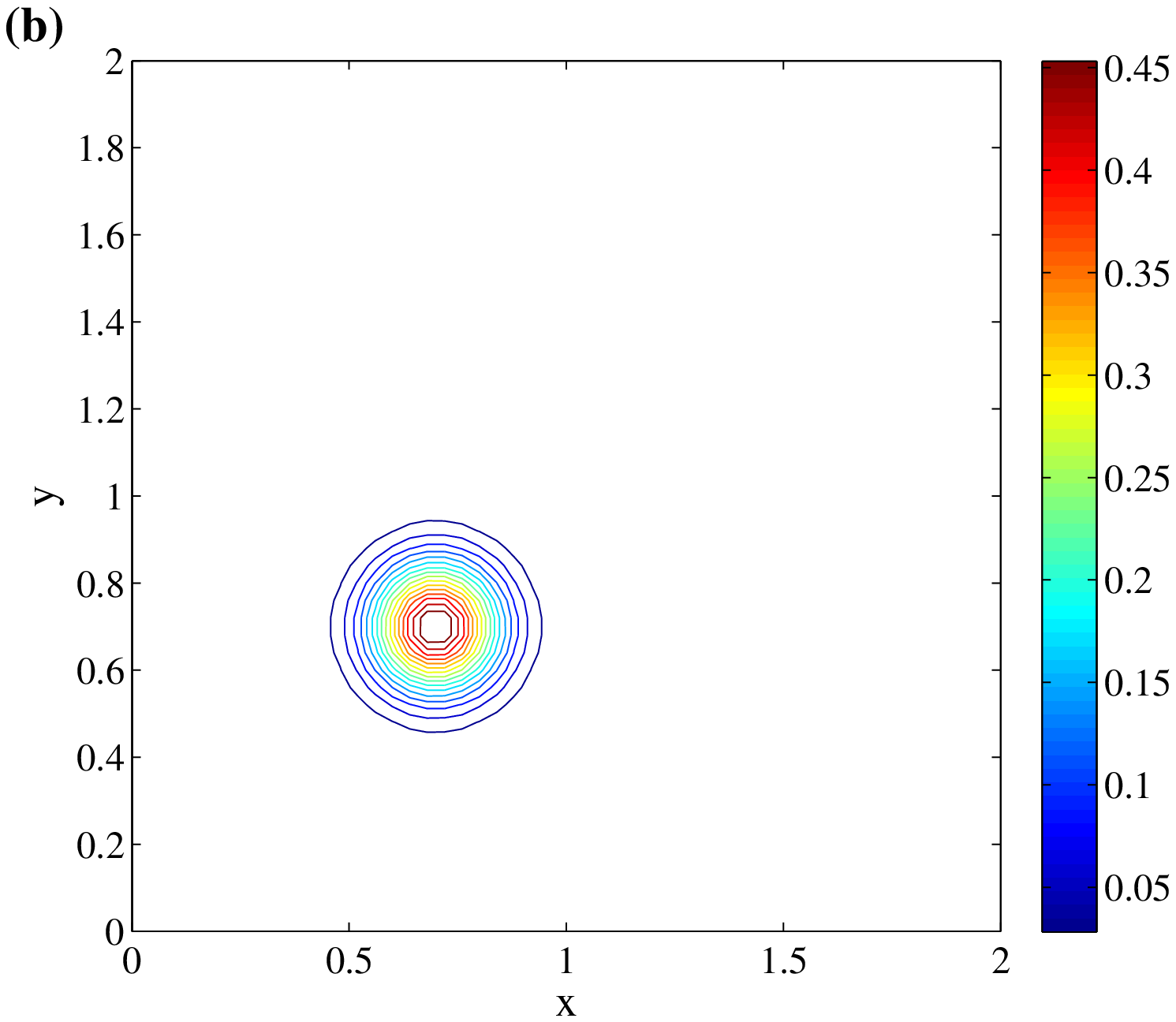}}\\
\subfigure{\includegraphics[width=2.7in]{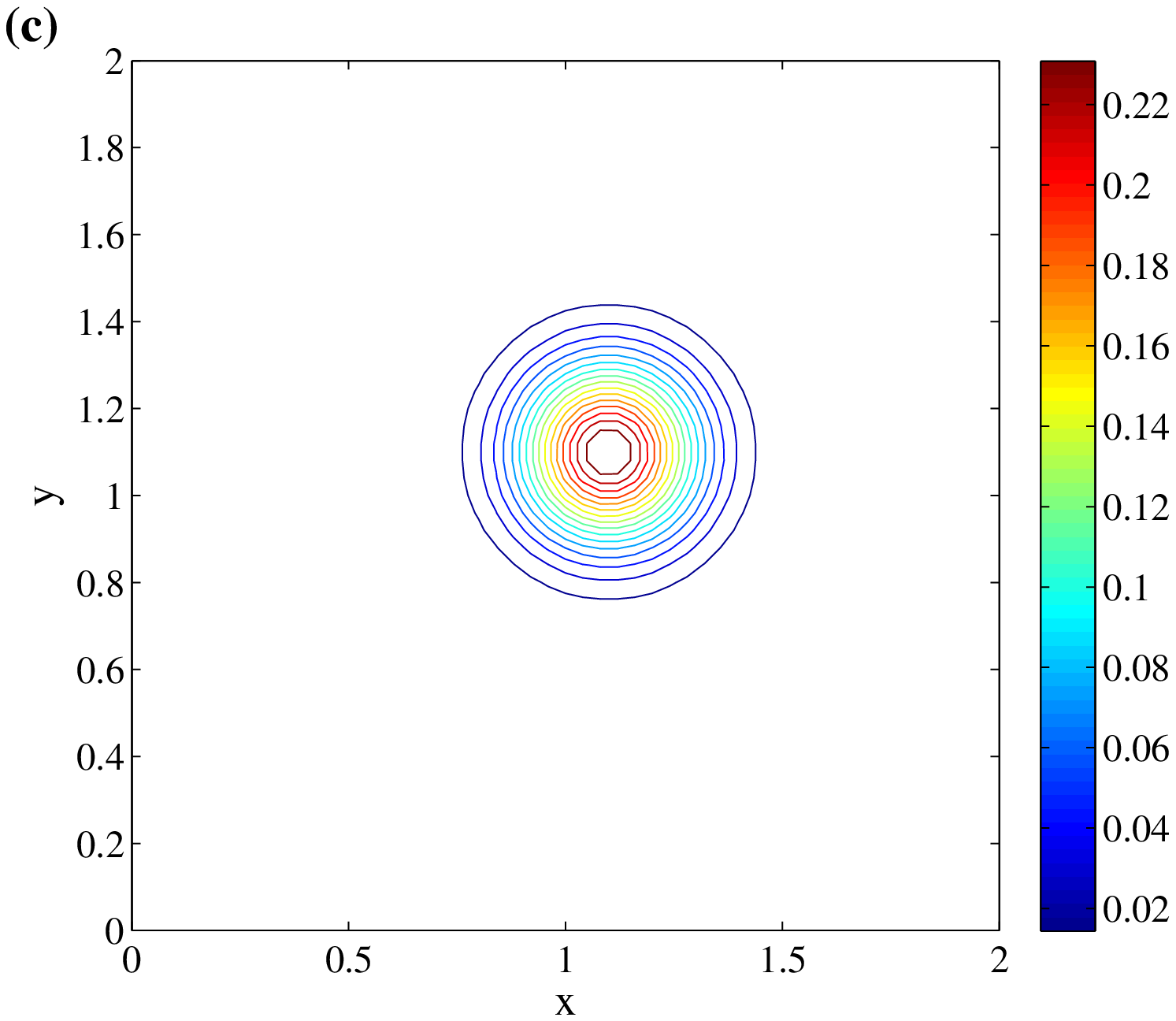}}\hspace{0.2in}
\subfigure{\includegraphics[width=2.7in]{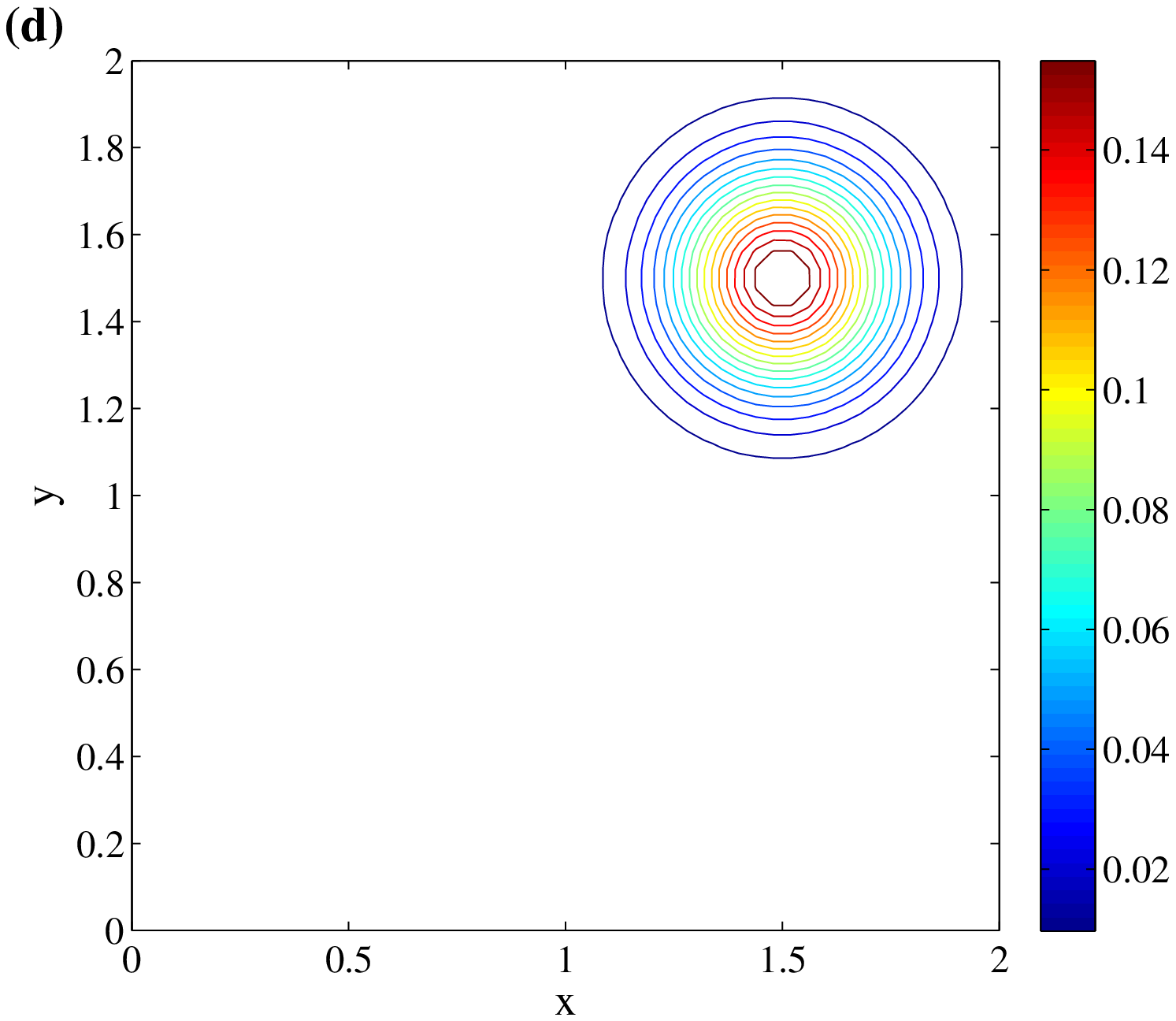}}
\caption{The contour plots of Gaussian pulse at $t=0$, $0.25$, $0.75$, $1.25$
 with $\tau=5.0\times10^{-3}$ and $M_x=M_y=50$.}\label{fig5}
\end{figure*}

\noindent
\emph{Example 6.6.} In this test, we simulate an unsteady propagation of a Gaussian pulse governed by
a classical 2D advection-dominated diffusion equation on a square domain $[0,2]\times[0,2]$ by using the scheme (\ref{ez09}),
which has been extensively studied \cite{ninep02,32,33,cirbf00}. The Gaussian pulse solution is expressed as
\begin{small}
\begin{equation*}
 u(x,y,t)\!=\!\frac{1}{1+4t}\exp\Bigg(-\frac{(x-\!\kappa_xt-\!0.5)^2}{\varepsilon_x(1+4t)}
 -\frac{(y-\!\kappa_yt-\!0.5)^2}{\varepsilon_y(1+4t)}\Bigg),
\end{equation*}
\end{small}
and the initial Gaussian pulse and boundary values are taken from the pulse solution.
Letting $\kappa_x=\kappa_y=0.8$, $\varepsilon_x=\varepsilon_y=0.01$, 
we display its true solution at $t=1.25$ with $M_x=M_y=50$ and the used lattice points on problem
domain in Fig. \ref{fig3}, which describe a pulse centred at $(1.5, 1.5)$ with a pulse height of $1/6$.
Using the same grid number together with $\tau=5.0\times10^{-3}$, we present the contour plots of the approximate solutions
at $t=0$, $0.25$, $0.75$, $1.25$ created by MCTB-DQM in Fig. \ref{fig5}. As the graph shows, the pulse
is initially centred at $(0.5, 0.5)$ with a pulse height of $1$, then it moves towards a position centred
at $(1.5, 1.5)$; during this process, its width and height appear to be continuously varying as the time goes by.
Besides, the last contour plot in Fig. \ref{fig5} coincides with the true solution plotted in Fig. \ref{fig3}.
Retaking $\tau=6.25\times10^{-3}$ and $M_x=M_y=80$, we compare our results with those obtained by some previous
algorithms as nine-point high-order compact (HOC) schemes \cite{ninep02,33}, Peaceman-Rachford ADI scheme (PR-ADI) \cite{pr08},
HOC-ADI scheme \cite{hoc03}, exponential HOC-ADI scheme (EHOC-ADI) \cite{ehoc05}, HOC boundary value method (HOC-BVM) \cite{bvm06},
compact integrated RBF ADI method (CIRBF-ADI) \cite{cirbf00}, coupled compact integrated RBF ADI method (CCIRBF-ADI) \cite{cirbf01},
and the Galerkin FEM combined with the method of characteristics (CGFEM) \cite{cgfem07}, at $t=1.25$ in Table \ref{tab5}. We implement
CGFEM on a quasi-uniform triangular mesh with the meshsize $2.5\times10^{-2}$ by using both Lagrangian P1 and P2 elements. Also,
average absolute errors are added as supplements to evaluate and compare their accuracy. As seen from Table \ref{tab5},
all of these methods are illustrated to be very accurate to capture the Gaussian pulse except the PR-ADI scheme;
besides, our method reaches a better accuracy than the others and even shows promise in treating
the advection-diffusion equations in the high P\'{e}clet number regime.\\

\begin{table*}
\centering
\caption{A comparison of global errors at $t=1.25$ with $\tau=6.25\times10^{-3}$ and $M_x=M_y=80$ for Example 6.6.}\label{tab5}
\begin{tabular}{llll}
\toprule Method    & Average Error & $e_2(M_x,M_y)$ &$e_\infty(M_x,M_y)$ \\
\midrule CIRBF-ADI \cite{cirbf00}   & 6.742e-06   & 2.197e-05    &1.703e-04    \\
         CCIRBF-ADI \cite{cirbf01}     & 5.989e-06   & 1.904e-05    &1.427e-04    \\
         Noye and Tan  \cite{33}        & 1.971e-05   & 1.280e-04    &6.509e-04    \\
         Kalita et al. \cite{ninep02}   & 1.597e-05   & 1.024e-04    &4.477e-04    \\
         PR-ADI  \cite{pr08}     & 3.109e-04  & 2.025e-03    &7.778e-03    \\
         HOC-ADI \cite{hoc03}    & 9.218e-06  & 5.931e-05    &2.500e-04    \\
         EHOC-ADI \cite{ehoc05}    & 9.663e-06    & 6.194e-05   &2.664e-04    \\
         HOC-BVM  \cite{bvm06}     & 9.4931e-06   & -   &2.4766e-04    \\
         CGFEM P1    & 6.3746e-05    & 1.8849e-04   &1.5731e-03    \\
         CGFEM P2    & 1.5667e-05    & 5.4061e-05   &5.8044e-04    \\
         MCTB-DQM    & 9.1512e-07    & 5.6996e-06   &2.2830e-05  \\
\bottomrule
\end{tabular}
\end{table*}

\noindent
\emph{Example 6.7.} In the last test, we consider the 2D space-fractional equations
(\ref{xz01})-(\ref{xz03}) on $[0,1]\times[0,1]$ with $\varepsilon_x=\varepsilon_y=1$, $\psi(x,y)=x^2(1-x)^2y^2(1-y)^2$,
and homogeneous boundary values. The source term is manufactured as
\begin{small}
\begin{align*}
    &f(x,y,t)=-e^{-t}x^2(1-x)^2y^2(1-y)^2\\
    &\ -\frac{2e^{-t}x^{2-\beta_1}y^2(1-y)^2}{\Gamma(3-\beta_1)}\Bigg(1
           -\frac{6x}{3-\beta_1}+\frac{12x^2}{(3-\beta_1)(4-\beta_1)}\Bigg)\\
    &\ -\frac{2e^{-t}x^2(1-x)^2 y^{2-\beta_2}}{\Gamma(3-\beta_2)}\Bigg(1
           -\frac{6y}{3-\beta_2}+\frac{12y^2}{(3-\beta_2)(4-\beta_2)}\Bigg)
\end{align*}
\end{small}
to enforce the analytic solution $u(x,y,t)=e^{-t}x^2(1-x)^2y^2(1-y)^2$. Letting $\beta_1=1.1$, $\beta_2=1.3$,
and $\tau=2.5\times10^{-4}$, we solve the problem via the FEM proposed by \cite{wx05}
and MCB-DQM, and compare their numerical results at $t=0.2$ in Table \ref{tab8}, where the P1 element and structured
meshes are adopted. The data indicate that DQ method converges towards the analytic solution as the grid numbers increase and
admits slightly better results than FEM. More importantly, the implemental CPU times of MCB-DQM are  less than
those of FEM, which confirm its  computing efficiency. \\ 

\section{Conclusion}\label{s6}
The ADEs are the subjects of active interest in mathematical physics and the related areas
of research. In this work, we have proposed an effective DQ method for such equations involving the
derivatives of fractional orders in time and space. Its weighted coefficients are calculated by making use of modified CTBs
and cubic B-splines as test functions. The stability of DQ method for the time-fractional ADEs in the context of $L_2$-norm is performed.
The theoretical condition required for the stable analysis is numerically  surveyed at length.
We test the codes on several benchmark problems and the outcomes have demonstrated that it outperforms some of the previously
reported algorithms such as BCFDM and FEM in term of overall accuracy and efficiency.

In a linear space, spanned by a set of proper basis functions as B-splines,
any function can be represented by a weighted combination of these basis functions.
While all basis functions are defined, the function remains unknown because the coefficients on the front of basis functions
are still unknown. However, when all basis functions satisfy Eqs. (\ref{xz17})-(\ref{xz18}), by virtue of linearity, it can be examined that the objective function
satisfies Eqs. (\ref{xz17})-(\ref{xz18}) as well. This is the essence of DQ methods, which guarantees their convergence.

Despite the error bounds are difficult to determine,
the numerical results illustrate that the spline-based DQ
method admits the convergent results for the fractional ADEs.
The presented approach can be generalized to the higher-dimensional and other
complex model problems arising in material science, structural and fluid mechanics, heat conduction, biomedicine,
differential dynamics, and so forth. High computing efficiency, low memory requirement, and the ease of
programming are its main advantages.

\begin{acknowledgements}
The authors are very grateful to the reviewers for their valuable comments and suggestions.
This research was supported by National Natural Science Foundations of China (Nos.11471262 and 11501450).
\end{acknowledgements}

\begin{table*}
\centering
\caption{A comparison of global errors at $t=0.2$ with $\tau=2.5\times10^{-4}$, $\beta_1=1.1$, and $\beta_2=1.3$ for Example 6.7.}\label{tab8}
\begin{tabular}{lcclclc}
\toprule
            Method  &  $M_x,\ M_y$   &$e_2(M_x,M_y)$ &Cov. rate &$ e_\infty(M_x,M_y)$ & Cov. rate & CPU time \\
\midrule FEM \cite{wx05} &10      &6.5781e-05  &-       &    1.5632e-04  &-      & 63.750(s)\\
                       & 15       &3.0082e-05  &1.9296  &    7.2433e-05  &1.8972 & 199.37(s)\\
                       & 20       &1.7376e-05  &1.9079  &    4.1321e-05  &1.9511 & 494.98(s)\\
                       & 25       &1.1404e-05  &1.8872  &    2.6599e-05  &1.9741 & 1062.7(s)\\
         MCB-DQM\!\!   & 10       &5.4217e-05  &-       &    1.4763e-04  &-      & 0.3726(s)\\
                       & 15       &2.6606e-05  &1.7556  &    6.9553e-05  &1.8562 & 1.1785(s)\\
                       & 20       &1.5559e-05  &1.8650  &    4.0088e-05  &1.9153 & 4.0598(s)\\
                       & 25       &1.0207e-05  &1.8892  &    2.6163e-05  &1.9124 & 11.267(s)\\
\bottomrule
\end{tabular}
\end{table*}

\begin{strip}
\noindent
\textbf{Appendix}:
\textbf{The explicit formulas of the fractional derivatives of cubic B-splines}

The fractional derivatives center at $x_{-1}$, $x_0$, and $x_1$:
\begin{align*}
{^{RL}_{x_0}}D^\beta_xB_{-1}(x)\!=\!\left\{ \begin{array}{l}
\frac{(1-\beta)(x-x_0)^{-\beta}}{\Gamma(2-\beta)}-\frac{3(x-x_0)^{1-\beta}}{\Gamma(2-\beta)h}
+\frac{6(x-x_0)^{2-\beta}}{\Gamma(3-\beta)h^2} -\frac{6(x-x_0)^{3-\beta}}{\Gamma(4-\beta)h^3}, \quad x \in [x_0,x_1)\\
\frac{(1-\beta)(x-x_0)^{-\beta}}{\Gamma(2-\beta)}-\frac{3(x-x_0)^{1-\beta}}{\Gamma(2-\beta)h}
+\frac{6(x-x_0)^{2-\beta}}{\Gamma(3-\beta)h^2} -\frac{6(x-x_0)^{3-\beta}}{\Gamma(4-\beta)h^3}
+\frac{6(x-x_1)^{3-\beta}}{\Gamma(4-\beta)h^3}, \quad x \in [x_1,x_M]
\end{array} \right.
\end{align*}

\begin{align*}
{^{RL}_{x_0}}D^\beta_xB_0(x)\!=\!\left\{ \begin{array}{l}
\frac{4(1-\beta)(x-x_0)^{-\beta}}{\Gamma(2-\beta)}-\frac{12(x-x_0)^{2-\beta}}{\Gamma(3-\beta)h^2}
+\frac{18(x-x_0)^{3-\beta}}{\Gamma(4-\beta)h^3}, \quad x \in [x_0,x_1)\\
\frac{4(1-\beta)(x-x_0)^{-\beta}}{\Gamma(2-\beta)}-\frac{12(x-x_0)^{2-\beta}}{\Gamma(3-\beta)h^2}
+\frac{18(x-x_0)^{3-\beta}}{\Gamma(4-\beta)h^3}-\frac{24(x-x_1)^{3-\beta}}{\Gamma(4-\beta)h^3}, \quad x \in [x_1,x_2)\\
\frac{4(1-\beta)(x-x_0)^{-\beta}}{\Gamma(2-\beta)}-\frac{12(x-x_0)^{2-\beta}}{\Gamma(3-\beta)h^2}
+\frac{18(x-x_0)^{3-\beta}}{\Gamma(4-\beta)h^3}-\frac{24(x-x_1)^{3-\beta}}{\Gamma(4-\beta)h^3}
+\frac{6(x-x_2)^{3-\beta}}{\Gamma(4-\beta)h^3}, \quad x \in [x_2,x_M]
\end{array} \right.
\end{align*}

\begin{align*}
{^{RL}_{x_0}}D^\beta_xB_1(x)\!=\!\left\{ \begin{array}{l}
\frac{(1-\beta)(x-x_0)^{-\beta}}{\Gamma(2-\beta)}+\frac{3(x-x_0)^{1-\beta}}{\Gamma(2-\beta)h}
+\frac{6(x-x_0)^{2-\beta}}{\Gamma(3-\beta)h^2} -\frac{18(x-x_0)^{3-\beta}}{\Gamma(4-\beta)h^3}, \quad   x \in [x_0,x_1)\\
\frac{(1-\beta)(x-x_0)^{-\beta}}{\Gamma(2-\beta)}+\frac{3(x-x_0)^{1-\beta}}{\Gamma(2-\beta)h}
+\frac{6(x-x_0)^{2-\beta}}{\Gamma(3-\beta)h^2} -\frac{18(x-x_0)^{3-\beta}}{\Gamma(4-\beta)h^3}
+\frac{36(x-x_1)^{3-\beta}}{\Gamma(4-\beta)h^3}, \quad x \in [x_1,x_2)\\
\frac{(1-\beta)(x-x_0)^{-\beta}}{\Gamma(2-\beta)}+\frac{3(x-x_0)^{1-\beta}}{\Gamma(2-\beta)h}
+\frac{6(x-x_0)^{2-\beta}}{\Gamma(3-\beta)h^2} -\frac{18(x-x_0)^{3-\beta}}{\Gamma(4-\beta)h^3}
+\frac{36(x-x_1)^{3-\beta}}{\Gamma(4-\beta)h^3}
-\frac{24(x-x_2)^{3-\beta}}{\Gamma(4-\beta)h^3},\quad  x \in [x_2,x_3)\\
\frac{(1-\beta)(x-x_0)^{-\beta}}{\Gamma(2-\beta)}+\frac{3(x-x_0)^{1-\beta}}{\Gamma(2-\beta)h}
+\frac{6(x-x_0)^{2-\beta}}{\Gamma(3-\beta)h^2} -\frac{18(x-x_0)^{3-\beta}}{\Gamma(4-\beta)h^3}
+\frac{36(x-x_1)^{3-\beta}}{\Gamma(4-\beta)h^3}-\frac{24(x-x_2)^{3-\beta}}{\Gamma(4-\beta)h^3}
+\frac{6(x-x_3)^{3-\beta}}{\Gamma(4-\beta)h^3}.\quad  x \in [x_3,x_M]
\end{array} \right.
\end{align*}

The fractional derivatives center at $x_{m}$ with $2\leq m\leq M+1$:

\begin{align*}
{^{RL}_{x_0}}D^\beta_xB_m(x)\!=\!\left\{ \begin{array}{l}
0,\qquad\qquad\qquad \, x \in [x_0,x_{m-2})\\
\frac{6(x-x_{m-2})^{3-\beta}}{\Gamma(4-\beta)h^3}, \quad x \in [{x_{m - 2}},{x_{m - 1}})\\
\frac{6(x-x_{m-2})^{3-\beta}}{\Gamma(4-\beta)h^3}-\frac{24(x-x_{m-1})^{3-\beta}}{\Gamma(4-\beta)h^3}, \quad x \in [{x_{m - 1}},{x_m})\\
\frac{6(x-x_{m-2})^{3-\beta}}{\Gamma(4-\beta)h^3}-\frac{24(x-x_{m-1})^{3-\beta}}{\Gamma(4-\beta)h^3}
+\frac{36(x-x_{m})^{3-\beta}}{\Gamma(4-\beta)h^3}, \quad x \in [{x_m},{x_{m + 1}})\\
\frac{6(x-x_{m-2})^{3-\beta}}{\Gamma(4-\beta)h^3}-\frac{24(x-x_{m-1})^{3-\beta}}{\Gamma(4-\beta)h^3}
+\frac{36(x-x_{m})^{3-\beta}}{\Gamma(4-\beta)h^3}-\frac{24(x-x_{m+1})^{3-\beta}}{\Gamma(4-\beta)h^3}, \quad x \in [{x_{m + 1}},{x_{m + 2}})\\
\frac{6(x-x_{m-2})^{3-\beta}}{\Gamma(4-\beta)h^3}-\frac{24(x-x_{m-1})^{3-\beta}}{\Gamma(4-\beta)h^3}
+\frac{36(x-x_{m})^{3-\beta}}{\Gamma(4-\beta)h^3}-\frac{24(x-x_{m+1})^{3-\beta}}{\Gamma(4-\beta)h^3}
+\frac{6(x-x_{m+2})^{3-\beta}}{\Gamma(4-\beta)h^3}, \quad x \in [x_{m+2},x_M]
\end{array} \right.
\end{align*}

which contain ${^{RL}_{x_0}}D^\beta_xB_{M-1}(x)$, ${^{RL}_{x_0}}D^\beta_xB_M(x)$, and ${^{RL}_{x_0}}D^\beta_xB_{M+1}(x)$ as special cases:

\begin{align*}
{^{RL}_{x_0}}D^\beta_xB_{M-1}(x)\!=\!\left\{ \begin{array}{l}
0,\qquad\qquad\qquad \, x \in [x_0,x_{M-3})\\
\frac{6(x-x_{M-3})^{3-\beta}}{\Gamma(4-\beta)h^3}, \quad x \in [x_{M - 3},x_{M - 2})\\
\frac{6(x-x_{M-3})^{3-\beta}}{\Gamma(4-\beta)h^3}-\frac{24(x-x_{M-2})^{3-\beta}}{\Gamma(4-\beta)h^3}, \quad x \in [{x_{M - 2}},x_{M-1})\\
\frac{6(x-x_{M-3})^{3-\beta}}{\Gamma(4-\beta)h^3}-\frac{24(x-x_{M-2})^{3-\beta}}{\Gamma(4-\beta)h^3}
+\frac{36(x-x_{M-1})^{3-\beta}}{\Gamma(4-\beta)h^3}, \quad x \in [x_{M-1},x_M]
\end{array} \right.
\end{align*}

\begin{align*}
{^{RL}_{x_0}}D^\beta_xB_M(x)\!=\!\left\{ \begin{array}{l}
0, \quad x \in [x_0,x_{M-2})\\
\frac{6(x-x_{M-2})^{3-\beta}}{\Gamma(4-\beta)h^3}, \quad x \in [x_{M-2},x_{M-1})\\
\frac{6(x-x_{M-2})^{3-\beta}}{\Gamma(4-\beta)h^3}-\frac{24(x-x_{M-1})^{3-\beta}}{\Gamma(4-\beta)h^3}, \quad x \in [x_{M-1},x_{M}]
\end{array} \right.
\end{align*}

\begin{align*}
{^{RL}_{x_0}}D^\beta_xB_{M+1}(x)\!=\!\left\{ \begin{array}{l}
0, \quad x \in [x_0,x_{M-1})\\
\frac{6(x-x_{M-1})^{3-\beta}}{\Gamma(4-\beta)h^3}. \quad x \in [x_{M-1},x_M]
\end{array} \right.
\end{align*}

\end{strip}

\end{document}